\documentclass[11pt]{article}

\usepackage[utf8]{inputenc}
\usepackage[scale=0.7]{geometry}
\usepackage[unicode]{hyperref}
\usepackage{parskip}
\usepackage{needspace}
\usepackage{enumerate}

\usepackage{newunicodechar}
\usepackage{amsmath}
\usepackage{amssymb}

\input{unicode-abbrv.sty}

\usepackage{tikz}
\usetikzlibrary{patterns}
\usetikzlibrary{decorations.markings}
\usepackage[margin=1cm]{caption}

\usepackage{bbm}

\usepackage{caption}
\usepackage[style=alphabetic,backend=bibtex]{biblatex}
\bibliography{Diplomarbeit}
\usepackage{booktabs}

\usepackage{makeidx}
\makeindex

\newcommand{\define}[1]{\emph{#1}\index{#1}}

\makeatletter
\def\qedsymbol{\ensuremath{\square}}
\def\qedEquationHack#1{\vskip-#1\hfill\qedsymbol\vskip\parskip}
\newenvironment{pr@@f}[1]{\textbf{\textit{Proof}#1.} }{}
\def\proof{\@ifnextchar[{\@with}{\@without}}
\def\@with[#1]{\pr@@f{ (#1)}}
\def\@without{\pr@@f{}}
\def\endproof{\hfill\qedsymbol\endpr@@f}
\def\proofx{\@ifnextchar[{\@with}{\@without}}
\def\@with[#1]{\pr@@f{ (#1)}}
\def\@without{\pr@@f{}}
\def\endproofx{\endpr@@f}
\makeatother

\newcommand{\Forward}{``$\Longrightarrow$''}

\newcommand{\Backward}{``$\Longleftarrow$''}

\makeatletter
\def\@begintheorem#1#2{\trivlist
   \item[\hskip \labelsep{\bfseries #1\ #2.}]\needspace{4\baselineskip}\mbox\newline\quote}
\def\@opargbegintheorem#1#2#3{\trivlist
      \item[\hskip \labelsep{\bfseries #1\ #2\ (#3).}]\needspace{4\baselineskip}\mbox\newline\quote}
\def\@endtheorem{\endquote\endtrivlist}
\makeatother

\newtheorem{theorem}{Theorem}[subsection]
\newtheorem{lemma}[theorem]{Lemma}
\newtheorem{proposition}[theorem]{Proposition}
\newtheorem{corollary}[theorem]{Corollary}

\newtheorem{definition}[theorem]{Definition}

\newenvironment{properties}{\begin{enumerate}[1)]}{\end{enumerate}}

\newcommand{\N}{\mathbb{N}}
\newcommand{\Z}{\mathbb{Z}}

\newcommand{\R}{\mathbb{R}}

\newcommand{\Union}{\bigcup}
\newcommand{\Continuous}{\mathcal{C}} % continuous functions, not complex numbers!
\newcommand{\supp}{\mathop{\text{supp}}}

\newcommand{\inverse}[1]{#1^{-1}}
\newcommand{\corona}{{}^*}

\newcommand{\defineiff}{{\ :\!\!\iff}}
\newcommand{\inject}{{\hookrightarrow}}
\newcommand{\true}{\text{true}}
\newcommand{\false}{\text{false}}

\newcommand{\1}{\mathbbm{1}}
\newcommand{\?}{\stackrel{?}}

\newcommand{\bN}{{\beta\N}}
\newcommand{\titlebN}{\texorpdfstring{$\bN$}{\unichar{"03B2}\unichar{"2115}}}
\newcommand{\kbN}{{\kappa(\bN)}}
\renewcommand{\d}{d^*}
\newcommand{\diagonalMap}{\Delta}

\newcommand{\aprich}{progression-rich} % rich in arithmetic progressions
\newcommand{\StoneCech}{Stone-\v{C}ech compactification}
\newcommand{\density}{upper Banach density}

\newcommand{\close}[1]{\overline{#1}}
\DeclareMathOperator{\closure}{cl}

\newcommand{\boundedSequences}{\ell^\infty(\N)}

\newcommand{\fold}[3]{%fold
{#1}_1 #2 {#1}_2 #2 \dots #2 {#1}_{#3}}
\renewcommand{\vector}[2]{\fold{#1},{#2}}
\newcommand{\partition}[1]{% set partition
{#1}_1 \uplus {#1}_2 \uplus \dots \uplus {#1}_r}
\newcommand{\sequence}[2]{({#1}_{#2})_{#2=1}^\infty}
\newcommand{\diagonal}[1]{(#1,#1,\dots,#1)}

\title{Formulating Szemerédi's Theorem in Terms of Ultrafilters}
\author{Heinrich-Gregor Zirnstein}
\date{March 24, 2011}

\begin{document}
\begin{titlepage}
\begin{center}
\large Universit\"at Leipzig\\
Fakult\"at f\"ur Mathematik und Informatik\\
Mathematisches Institut
%\\[2.5cm]
\\[5cm]
\large Diplomarbeit \\[1.5cm]
\LARGE Formulating Szemerédi's Theorem\\ in Terms of Ultrafilters
\\[7cm]
%\\[2.8cm]
%\large Betreuender Hochschullehrer:\\
%Professor Dr. Andreas Thom\\[6cm]
\end{center}
\normalsize
\begin{tabular}{p{6.8cm}ll}
Leipzig, den 22. April 2012 & Vorgelegt von: & Heinrich-Gregor Zirnstein\\
& Studiengang: & Mathematik, Diplom \\
& Betreuer: & Professor Dr. Andreas Thom
\end{tabular}
\end{titlepage}

\newpage
\tableofcontents
\thispagestyle{empty}
\newpage

%%%%%%%%%%%%%%%%%%%%%%%%%%%%%%%%%%%%%%%%%%%%%%%%%%%%%%%%%%%%%%%%%%%%%%
\section*{Preface}
\def\myskip{\vskip0pt}
\myskip
\addcontentsline{toc}{section}{Preface}
%%%%%%%%%%%%%%%%%%%%%%%%%%%%%%%%%%%
\subsubsection*{Introduction}
\myskip
A famous theorem by van der Waerden \parencite{vanderWaerden:1927} asserts that if you color the natural numbers with, say, five different colors, then you can always find arbitrarily long sequences of numbers that have the same color and that form an arithmetic progression; the same is true for any other count of colors.

Since coloring the natural numbers is the same as partitioning them, $\N=\partition C$, one also says that arithmetic progressions are \emph{partition regular}: no matter how you divide the natural numbers, one of the parts will contain arithmetic progressions of arbitrary length. The branch of mathematics concerned with the study of partition regularity is called \emph{Ramsey theory}.

Van der Waerden's theorem appears to be a beautiful statement from Ramsey theory that only requires elementary mathematics to be understood and proven. While this is indeed the case, there is a different point of view based on \emph{ultrafilters}, strange objects from point-set topology whose very existence is linked to the axiom of choice and hence somewhat mysterious.

It turns out that partition regular properties correspond to ultrafilters with special properties. Furthermore, the existence of an ultrafilter corresponding to van der Waerden's theorem, and hence the theorem itself, can be proven by exploiting a semigroup structure on the space of ultrafilters.

A remarkable generalization of van der Waerden's theorem is \emph{Szemerédi's theorem} \parencite{Szemeredi:1975p1591}, \parencite{Green:2007}, which is much deeper and asserts that every subset of natural numbers with positive density contains arithmetic progressions of arbitrary length. The density of a set $A\subseteq\N$ is defined as the limes superior of the ratios $|A\cap[1,N]|/N$ as $N \to \infty$. It is easy to see that in any partition $\N=\partition C$, one of the parts must have positive density, so this is indeed a generalization.

The aim of this diploma thesis is to give an interpretation of Szemerédi's theorem in terms of ultrafilters as well. Namely, while van der Waerden's theorem is equivalent to the existence of a single ultrafilter with special properties, we will show that Szemerédi's theorem is equivalent to the existence of not just one, but of many such ultrafilters. In fact, we will define a measure on the space of ultrafilters $\bN$ and deduce that, with respect to this measure, \emph{almost all} ultrafilters must have the special properties needed for van der Waerden's theorem.

This novel interpretation does not necessarily yield a proof of Szemerédi's theorem, however. In particular, we will show that the ultrafilter proof of van der Waerden's theorem is far too weak to imply Szemerédi's theorem. The reason is that the set of special ultrafilters exhibited by the argument has measure zero, even though it is infinite.

%%%%%%%%%%%%%%%%%%%%%%%%%%%%%%%%%%%
\subsubsection*{About this text}
\myskip
When writing this thesis, I have tried to keep the material as accessible as possible; everything is explained and motivated thoroughly. In proofs, I have strived for both clarity and detail, although that makes them somewhat lengthy on paper. They can be skipped on first reading.

The text is largely self-contained, the only prerequisite is familiarity with basic point-set topology and measure theory, for example as presented in \textcite{Janich:2005p1530} and the first chapters of \textcite{Elstrodt:2005p567}. At one point, however, we will make use of the Riesz representation theorem, which is discussed in the more advanced parts of the mentioned book by Elstrodt.

%%%%%%%%%%%%%%%%%%%%%%%%%%%%%%%%%%%
\subsubsection*{Chapter overview}
\myskip
Here a synopsis of the individual chapters of this diploma thesis.

\begin{description}
\item[Chapter \ref{chapter-ultrafilters}] introduces the space of ultrafilters $\bN$ and its topological properties. We also explain limits along ultrafilters.

\item[Chapter \ref{chapter-ramsey}] recalls the relevant notions from Ramsey theory and gives a proof that all partition regular properties correspond to ultrafilters.

\item[Chapter \ref{chapter-algebra}] defines the addition of ultrafilters, which turns the space of ultrafilters $\bN$ into a left topological semigroup. The notation $A-p$ for ultrafilter shifts is also introduced. We will study idempotent ultrafilters $p=p+p$ and prove Hindman's theorem about IP-sets. Finally, we will collect information about the ideals of the  semigroup $\bN$ to the point that we can give a proof of van der Waerden's theorem.

\item[Chapter \ref{chapter-measure}] introduces a family of counting measures on the space of ultrafilters $\bN$. As an application, we will reproduce Beiglböck's proof \parencite{Beiglbock:2009p581} of Jin's theorem about the size of sets of differences $A-B$. Finally, we give and prove the interpretation of Szemerédi's theorem in terms of counting measures and ultrafilters. Also, we will argue that the proof of van der Waerden's theorem from the previous chapter cannot imply Szemerédi's theorem.
\end{description}

The material for the first three chapters is taken mainly from \textcite{Hindman:1998p531} and \textcite{Bergelson:2003p505}.
\myskip
%%%%%%%%%%%%%%%%%%%%%%%%%%%%%%%%%%%
\subsubsection*{Acknowledgments}
\myskip
I thank my supervisor Andreas Thom for inspiring discussions and for nurturing an exceptional environment for doing mathematics.

%%%%%%%%%%%%%%%%%%%%%%%%%%%%%%%%%%%%%%%%%%%%%%%%%%%%%%%%%%%%%%%%%%%%%%
\newpage
\section{The space of ultrafilters \titlebN}
\label{chapter-ultrafilters}
%%%%%%%%%%%%%%%%%%%%%%%%%%%%%%%%%%%
\subsection{Missing natural numbers?}

We begin with the definition of ultrafilters.

\begin{definition}[Ultrafilter]
Let $X$ be any set. A nonempty collection $p$ of subsets of $X$, $p\subseteq 2^X$, is called an \define{ultrafilter} on $X$ if it satisfies the following properties:
\begin{properties}
	\item $p \not\owns \emptyset$.
	\item If $p \owns A$ and $A \subseteq B$, then $p \owns B$.
	\item If $p \owns A$ and $p \owns B$, then $p \owns A\cap B$.
	\item Either $p \owns A$ or $p \owns A^c$ for all sets $A$. (We write $A^c = X\setminus A$ for the complement.)
\end{properties}
We denote the \emph{set of ultrafilters}\index{ultrafilter!set of ultrafilters@set of ultrafilters $\beta X$} on $X$ with $\beta X$\footnote{Actually, $\beta X$ is the standard notation for the \StoneCech, which coincides with the set of ultrafilters if $X$ is a discrete topological space. See Definition \ref{definition-stone-cech}.}.
\end{definition}

We are mainly interested in ultrafilters over the natural numbers, so we usually take $X=\N$. Occasionally, we will be interested in other discrete sets like $X = \N\times\N$ as well.

The intuition behind ultrafilters is that they behave like ``the missing elements'' of the set $X$. Well, it is probably news to you there might be anything missing from $X$, but if you look at the four properties above and flip the membership symbol from ``$p \owns$'' to ``$p \in$'', you will suddenly notice that these properties are a reasonable axiomatization of set membership. For instance, the third property would read
\[ \text{If $p \in A$ and $p \in B$, then $p \in A\cap B$}, \]
which is just the definition of set intersection. In a sense, we are adding new elements, or ``points'', to the set $X$. These new points $p$ are specified by the collection of sets $p = \{A,B,\dots\}$ in which they are ``contained'' in. We will make this intuition rigorous in Section \ref{section-topology-bN}, when we study the \emph{space of ultrafilters} $\beta X$.

In this light, our first examples of ultrafilters are the so called \emph{principal ultrafilters}\index{ultrafilter!principal}
\[ \tilde x = \{ A\subseteq X : x\in A \}, \]
which simply correspond to the original points of $X$. Ultrafilters that are not of this form are the interesting ``new'' points, they are called \emph{non-principal ultrafilters}\index{ultrafilter!non-principal}.

Of course, the question is whether there exist any non-principal ultrafilters at all. The answer is ``yes''; we will construct non-principal ultrafilters shortly using Zorn's lemma. But note that their existence is independent of the ZF axioms of set theory, so we really need to use the Axiom of Choice here; see also \textcite{Schechter:2008p1592}. This has the unfortunate consequence that we cannot enumerate the collection of sets defining an ultrafilter in a meaningful, algorithmic way. Hence, ultrafilters will remain somewhat mysterious objects.

%%%%

But before constructing ultrafilters, let us define the related notion of \emph{filter}, which corresponds not to a single point, but to a set of points.

\newcommand{\F}{\mathcal F}
\begin{definition}[Filter]
Let $X$ be a set. A \define{filter} $\F$ on $X$ is a nonempty collection of subsets of $X$ that fulfills the properties 1--3\ from the definition of ultrafilters, but not necessarily the condition 4.
\end{definition}

This time, imagine replacing ``$\F\owns$'' with ``$\F\subseteq$'' and observe that the conditions 1--3 are a reasonable axiomatization of the notion of subset. In this light, condition 4 singles out the subsets with just one element, i.e.\ the points.

Let us give some examples of filters. Like before, the \emph{principal filters}\index{filter!principal} are the ones that come from ordinary subsets $F$ of the set $X$:
\[  \tilde F = \{ A\subseteq X : F \subseteq A\}  .\]
Similarly, the filters that are not of this form are called \emph{non-principal filters}\index{filter!non-principal}. The simplest example is the \emph{Fréchet filter}, which consists of those sets whose complement is finite
\[  \F = \{ A\subseteq X : A^c \text{ is a finite set}\} ,\]
or in the case of $X=\N$
\[  \F = \{ A\subseteq \N : \{n,n+1,\dots\} \subseteq A \text{ for some } n\in\N \}  .\]
This is indeed a filter, and non-principal one because the intersection of all of its sets is empty, which cannot happen for principal filters.

%%%%

Back to ultrafilter construction. In a sense, condition 4 from the ultrafilter definition is also a condition of maximality. Clearly, a filter $\F$ cannot contain both a set $A$ and its complement $A^c$, because then we would have $\F \owns A\cap A^c = \emptyset$, in violation of property 1. But it may well be that a filter contains neither. In contrast, ultrafilters must be large (``ultra'') and always contain one of them. We will show that any filter can be extended to an ultrafilter.

\begin{lemma}[Extending filters to ultrafilters]
\label{extend-filter}
Let $\F$ be a filter on the set $X$. Then, there exists an ultrafilter $p$ that contains all the sets from $\F$ and many more, i.e.\ $\F \subseteq p$.
\end{lemma}

In our interpretation, this lemma says that ``every set $\F$ contains a point $p$''. Applying it to the Fréchet filter shows the existence of non-principal ultrafilters.

\begin{corollary}[Existence of non-principal ultrafilters]
There exist non-principal ultrafilters $p$. Moreover, every non-principal ultrafilter extends the Fréchet filter.
\end{corollary}
\begin{proof}
As the previous lemma says, there exists an ultrafilter $p$ extending the Fréchet filter. Since the intersection of all sets $A\in p$ is empty, it cannot be principal.

Now, note that only principal ultrafilter $p$ may contain finite sets. Hence, given any set $A$ whose complement $A^c$ is finite, a non-principal ultrafilter $p$ must contain the set $A$ but not the set $A^c$. This means that $p$ extends the Fréchet filter.
\end{proof}

For reasons of economy, we now prove a generalization of the lemma above that involves a predicate on sets $\phi$. Setting $\phi(A) = \text{``$A$ is nonempty''}$ will recover the original statement.

%:	ultrafilter construction lemma

\begin{lemma}[Ultrafilter construction]
\index{ultrafilter construction lemma}
\label{prop-ultrafilter-construction}
Let $X$ be a set and $\F$ be a filter on $X$. Furthermore, let $\phi : 2^X \to \{\true,\false\}$ be a predicate on subsets of $X$ that has the following properties.
\begin{properties}
\item $\phi(A)$ for all $A\in\F$.
\item If $\phi(A)$ and $A\subseteq B$, then $\phi(B)$.
\item If $\phi(A)$ and $A=A_1\uplus A_2$ a disjoint union, then $\phi(A_1)$ or $\phi(A_2)$ or both.
\end{properties}
In other words, the predicate cannot be ``destroyed'' by partitioning a set\footnote{The set $\{A\subseteq X : \phi(A)\}$ is sometimes called a \define{superfilter}.}.

Then, there exists an ultrafilter $p$ with $\F \subseteq p$ and $\phi(A)$ for all $A\in p$.
\end{lemma}
\begin{proof}
Consider all filters $\mathcal G \supseteq \F$ whose member sets $A\in \mathcal G$ all have the property $\phi(A)$. By Zorn's lemma, there exists a filter $p$ among these which is maximal with respect to inclusion. We want to show that this is an ultrafilter.

Assume that there were a set $A$ such that neither $A$ itself nor its complement $A^c$ are members of $p$. We show that in this case, it is possible to extend the filter $p$ by one of these sets, in contradiction to maximality.

Let $p[A]$ denote the filter ``generated'' by the filter $p$ and the set $A$, i.e.\ the collection
\[  p[A] = p \cup \{S\subseteq X: \text{ there exists } B\in p \text{ with } A\cap B \subseteq S\}\]
Likewise for $p[A^c]$. It is easy to check that these collections are indeed filters. In particular, since neither the set $A$ nor its complement $A^c$ are members of the filter $p$, none of the intersections $A\cap B$ and $A^c\cap B$ with $B\in p$ can be empty. Hence, neither collection contains the empty set.

We have to prove that all members sets of at least one of these filters satisfy the predicate $\phi$. Thanks to requirement 2, it is enough to show that
\[ \text{either}\quad \forall B\in p.\ \phi(A\cap B) \quad\text{or}\quad \forall C\in p.\ \phi(A^c\cap C) \quad\text{or both}.\]
In words: the set $A$ or the set $A^c$ should only have ``good'' intersections with all the members of the filter $p$.

%: FIGURE: ultrafilter construction lemma
\begin{figure}
\begin{center}
\includegraphics{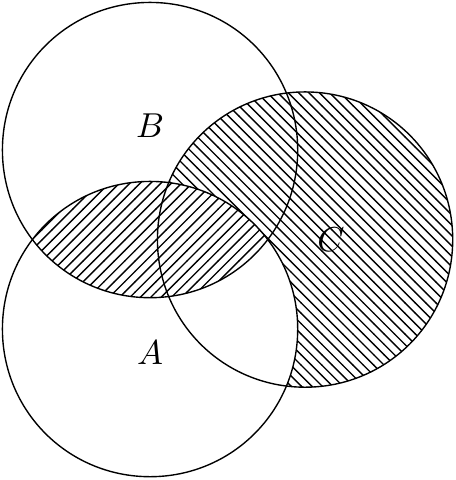}
\end{center}
\caption{The predicate $\phi$ being false on the shaded regions $A\cap B$ and $A^c \cap C$ implies that it is also false on the intersection $B\cap C$.
\label{pic-ultrafilter-construction}}
\end{figure}

Assume that this is not true, i.e.\ that there are counterexamples $B$ and $C$ such that the statements $\phi(A\cap B)$ and $\phi(A^c\cap C)$ are false, see Figure \ref{pic-ultrafilter-construction}. Due to requirement 2 again, the statements $\phi(A\cap B\cap C)$ and $\phi(A^c\cap C\cap B)$ must also be false. But the contrapositive of requirement 3 tells us that
\[  \phi((A\cap B\cap C) \uplus (A^c\cap C\cap B)) = \phi(B\cap C) = \false \]
as well, in contradiction to $B\cap C\in p$ and $\phi(B \cap C) = \true$.

In other words, at least one of the filters $p[A]$ or $p[A^c]$ extends the filter $p$, even though the latter was supposed to be maximal.
\end{proof}

In the end, we want to use ultrafilters on the natural numbers to prove statements about the natural numbers themselves. The passage from ultrafilters back to natural numbers can be captured in the following simple observation, whose value lies in the mental image it provokes: If a set $A$ contains a ``non-principal point'' $p \owns A$, then it also contains an ordinary point $n \in A$, $n \in \N$.
\begin{proposition}[Permanence principle]
\index{permanence principle}
\label{prop-permanence-principle}
Let $A$ be a set of natural numbers. For any ultrafilter $p$ on the natural numbers $\N$, we have
\[ A \in p \implies A \owns n \text{ for some natural number } n\in\N \]
Moreover, if the ultrafilter $p$ is non-principal, then this can be strenghened to
\[ A \in p \implies \begin{array}{l}\text{For any number $m\in\N$, there exists a number $n\in\N$}\\[-5pt]\text{such that $n>m$ and } A \owns n.\end{array}\]
\end{proposition}
\begin{proof}
Ultrafilters do not contain the empty set. Non-principal ultrafilters are contained in the Fréchet filter.
\end{proof}

Let us close this section with a lemma that we will use in the next section to make the intuition of ultrafilters as points rigorous.

\begin{proposition}[Ultrafilters behave like points]
\label{prop-ultrafilters-points}
For any ultrafilter $p$ on some set $X$ and any sets $A,B\subseteq X$, we have the following equivalences:
\begin{enumerate}
\item $p\owns A^c$ $\iff$ not $p\owns A$ 
\item $p\owns A\cap B$ $\iff$ $p\owns A$ and $p\owns B$ 
\item $p\owns A\cup B$ $\iff$ $p\owns A$ or $p\owns B$
\end{enumerate}
\end{proposition}
\begin{proofx}
\begin{enumerate}
\item By definition, one of $A$ and $A^c$ must be a member of $p$, but since $A\cap A^c=\emptyset\not\in p$, only one of them can be.
\item The implication from right to left holds by definition. The other direction follows trivially from $A\cap B\subseteq A$ and $A\cap B\subseteq B$.
\item De Morgans' law reduces this to the first two cases:
\begin{align*}
p\owns A\cup B &\iff p\owns (A^c\cap B^c)^c
\\ &\iff \text{not $(p\owns A^c$ and $p \owns B^c)$} \iff p\owns A \text{ or } p\owns B.
\end{align*}\qedEquationHack{0.98cm}
\end{enumerate}
\end{proofx}

%%%%%%%%%%%%%%%%%%%%%%%%%%%%%%%%%%%
\subsection{Topology of \titlebN}
\label{section-topology-bN}
We now want to make rigorous our intuition of ultrafilters as points and of filters as sets of points.

For simplicity, the following discussion will be about ultrafilters on the set of natural numbers $\N$. It also applies to $\Z$, $\N\times\N$ and other discrete sets.\footnote{When the set carries a non-trivial topology that needs to be taken into account, the \StoneCech\ (Definition \ref{definition-stone-cech}) is more appropriate.}

To interpret ultrafilters as \emph{points}, we can simply decree that they are the points of the \emph{space of ultrafilters} $\bN$. The ordinary natural numbers $\N$ can be viewed as points of $\bN$ as well, by thinking of them as principal ultrafilters. In other words, there is an inclusion $\N\inject\bN$.

But now, there are two possible notions of \emph{sets of points}: on one hand, we have arbitrary subsets of $\bN$, in the sense of set theoretical collections of points. On the other hand, we want to think of \emph{the filters} as being the subsets of the space $\bN$. The solution to this dilemma is to introduce a \emph{topology} on the space $\bN$. Arbitrary sets will remain arbitrary, but the filters will correspond to topologically meaningful sets; we will see that they correspond to the \emph{closed sets} of the space $\bN$.

Furthermore, it will turn out that the space $\bN$ is a \emph{compact} topological space, unlike the natural numbers $\N$ themselves. So, in a sense, there \emph{were} points missing from $\N$, namely those that correspond to limit points, without whom $\N$ cannot be compact.

At the end of this section, we will give short summary of the topology and illustrate it with a picture.

%:	open sets

To define the topology of the space $\bN$, we first define the closure of a set of natural numbers.
\begin{definition}[Closure of a set of natural numbers in $\bN$]
\label{definition-closure}
Let $A\subseteq\N$ be a set of natural numbers. Its \emph{closure $\close A$ in $\bN$}\index{closure!of a set of natural numbers} is defined to be the set of all ultrafilters $p$ that it belongs to:
\[  p\in\close A \defineiff p\owns A  .\]
\end{definition}

In other words, our intuition that the ultrafilter $p$ is a ``point of $A$'' is now manifest in the fact that $p$ is a point of the closure $\close A$.

Note that we reserve the notation ``$\close A$'' for the closure of sets of natural numbers $A \subseteq \N$. The topological closure of arbitrary sets $X \subseteq \bN$ will be denoted with ``$\closure(X)$''\index{closure!topological closure}. We will justify that $\close A = \closure(A)$ shortly.

In a very confusing move, we now define the closures of sets of natural numbers to be the basic \emph{open} sets of our topology.

\begin{definition}[Topology of $\bN$]
\def\a{\alpha}
\def\A{\mathcal A}
A set of the form $\close A$ is called a \define{basic open set}. A set of ultrafilters $U\subseteq\bN$ is called \emph{open} if and only if it is a union of basic open sets, 
\[  U = \bigcup_{\a\in\A} \close A_\a \text{ for some collection of sets } A_\a \subseteq \N  .\]
\end{definition}

We just defined these \emph{closures} to be \emph{open}; but are they not supposed to be \emph{closed}? The solution to this conundrum is that they are \emph{both}, they are both closed and open, they are \emph{clopen}. This might be unfamiliar to you if you are used to the topology of the real numbers where the open and the closed sets are quite distinct; but even there, two special sets are clopen, namely the whole space and the empty set. Well, here we have a topological space where many interesting sets are clopen.

The following proposition justifies that sets of the form $\close A$ are both open and closed, and that our definition of ``open set'' adheres to the axioms of topology.

\begin{proposition}[Closure commutes with boolean algebra]
\label{prop-boolean-algebra}
For any sets $A,B\subseteq\N$, the following equalities hold in $\bN$:
\begin{enumerate}
\item $\close {A^c} = {\close A}^c$
\item $\close {A\cap B} = \close A \cap \close B$
\item $\close {A\cup B} = \close A \cup \close B$
\end{enumerate}
\end{proposition}
\begin{proofx}
All these equalities follow directly from Proposition \ref{prop-ultrafilters-points} ``Ultrafilters behave like points'', which has been stated in such a way that we now only have to flip the ``$\in$'' symbol forth and back. Here the reasoning for the first equality; the others are entirely similar:
\[ p\in \close {A^c} \iff p\owns A^c \iff \text{not }(p\owns A) \iff \text{not }(p\in\close A) \iff p\in{\close A}^c . \]
\qedEquationHack{0.98cm}
\end{proofx}

We still have to justify that the ``closure'' $\close A$ is indeed equal to the smallest closed set $\closure(A)$ of ultrafilters that contains the set $A \subseteq \N$.
\begin{lemma}[Closure is topological closure]
For every set of natural numbers $A\subseteq\N$, we have $\close A = \closure(A)$.
\end{lemma}
\begin{proof}
Without loss of generality, we can write the closure as an intersection of basic open (closed) sets $\closure(A) = \bigcap_{\alpha\in\mathcal A} \close B_\alpha$. But the ultrafilter axioms show that $A \subseteq \close B_\alpha$ implies $\close A \subseteq \close B_\alpha$ and the intersection must be equal to its smallest member $\close A$.
\end{proof}

%:	filters = closed sets

Having defined the topology on $\bN$, we can now state the correspondence between filters and closed sets.

\begin{proposition}[Correspondence of closed sets and filters]
\label{prop-filter-closed}
The nonempty closed sets $F\subseteq\bN$ are in one-to-one correspondence with the filters $\F$ on $\N$. Moreover, the subset relation $F\subseteq G$ corresponds to filter inclusion $\F \supseteq \mathcal G$.

In particular, a filter $\F$ corresponds to the intersection $F$ of all the sets that we used to define it as a collection,
\[  F = \bigcap_{A\in\F} \close A  .\]
\end{proposition}

\newcommand\filter{f}
\newcommand\closed{c}
\begin{proof}
With this intuition in mind, let us introduce two maps $\closed$ and $\filter$ that map a filter to its closed set and vice-versa:
\begin{align*}
\closed(\F)  &= \bigcap_{A\in\F} \close A\\
\filter(F)   &= \{ A\subseteq\N : F \subseteq \close A\}
\end{align*}
Since arbitrary intersections of closed sets are closed, the set $\closed(\F)$ is indeed a closed set. Likewise, since closure commutes with boolean algebra, $\filter(F)$ is indeed a filter. It is also obvious that these two maps interchange subset relation and filter inclusion.

Note that every ultrafilter $p$ is a member of the set $\closed(p)$; this is obvious from the definition of closure. Since every filter $\F$ can be extended to an ultrafilter $p\supseteq \F$, we have $p\in\closed(p)\subseteq c(\F)$, so that the closed set $c(\F)$ is, in fact, \emph{nonempty}.

It remains to be shown that these maps are \emph{inverse} to each other.

Let us consider the case $F \?= \closed(\filter(F))$ first. From the definition of set intersection, it is obvious that $F \subseteq \closed(\filter(F))$. But since the closures $\close A$ form a \emph{basis of our topology}, every closed set $F$ can be written as an intersection $F=\bigcap_{A\in\mathcal A}\close A=\closed(\mathcal A)$ for some collection $\mathcal A$ (which, however, is not necessarily a filter). This means $\mathcal A \subseteq \filter(F)$, which implies $F =\closed(\mathcal A) \supseteq \closed(\filter(F))$, and we are done.

Now the case $\F \?= \filter(\closed(\F))$. Once again, it is obvious that $\F \subseteq \filter(\closed(\F))$. If we can show that $\closed$ is \emph{injective} in the sense that $\F \subset \F'$ for some other filter $\F'$ implies $\closed(\F) \not= \closed(\F')$, then we can apply the previous case to the closed set $F=\closed(\F)$ and conclude from injectivity that
\[  F=\closed(\filter(F)) \implies \closed(\F)=\closed(\filter(\closed(\F))) \implies \F = \filter(\closed(\F)) .\]

Alright then, let us prove that the map $\closed$ is injective. Intuitively, we have to show that $\closed(\F')$ is a strictly smaller set than $\closed(\F)$, so we have to construct an ultrafilter $p$ which is a member of the latter, but not of the former closed set. Since $\F'$ is a larger filter than $\F$, we can find a set $B$ with $B\in\F'$ but $B\not\in\F$. Now, the idea is that the closed set $\closed(\F)$ has nonempty intersection with the complement $\close B^c$ whereas the closed set $\closed(\F')$ must be fully contained in the closure $\close B$. To show that, and to conclude the proof, all we have to do is to construct an ultrafilter $p\in\closed(\F)\cap\close B^c$.

Of course, to construct the ultrafilter $p$, we apply the \emph{ultrafilter construction lemma} (\ref{prop-ultrafilter-construction}) to the filter $\F$, using the predicate ``$\phi(A) = A\cap B^c \text{ is nonempty}$''. As we will see in a moment, this will have the effect that $p$ contains both the filter $\F$ and the set $B^c$. The most important condition to verify is that this predicate is true for all sets $A\in\F$. But if that were not the case, that is if one set $A\in\F$ fulfilled ``$A\cap B^c \text{ is empty}$'', then we would have $A\subseteq B$ which would imply $B\in\F$ in contradiction to our choice of $B$. The other conditions are easy.

Now, all that remains is to check that indeed $p\in\closed(\F)\cap\close B^c$. Since $p$ is an extension of $\F$, we know that $p\in c(p)\subseteq c(\F)$. Also, $p$ must be contained in either $\close B$ or $\close B^c$, but since $\phi(B)=\false$, we have $p\in\close B^c$ as desired.
\end{proof}

%:	compactness

The fact that the closed sets corresponding to filters are \emph{nonempty} is closely related to the fact that our space $\bN$ is compact.

\begin{proposition}[Compactness of $\bN$]
The space of ultrafilters $\bN$ is a compact topological space.
\end{proposition}
Remember that compactness means that all open covers have finite subcovers and that the space is Hausdorff. To prove the proposition, let us reformulate the first criterion.
\begin{definition}[Finite intersection property]
\def\a{\alpha}
\def\A{\mathcal A}
A topological space $X$ has the \define{finite intersection property} if, for every arbitrary collection of closed sets $\{C_\a\subseteq X : \a\in\A\}$, the following equivalence holds:
\begin{align*}
\bigcap_{\a\in\A} C_\a \text{ is empty } \iff
\begin{array}{l}
\text{ there exists a finite set } I\subseteq\A
\\[-5pt] \text{ such that } \bigcap_{\a\in I} C_\a \text{ is empty } .
\end{array}
\end{align*}
\end{definition}

Take complements to see that the finite intersection property is equivalent to the usual property concerning open covers.

\begin{proof}[Compactness of $\bN$]
\def\a{\alpha}
\def\A{\mathcal A}
First, we show that the space $\bN$ is Hausdorff. To see this, consider two ultrafilters $p \neq q$. Since they are different, there must be a set $A$ such that $A \in p$ but $A \not\in q$. This means $p \in \close A$ and $q \in \close{A}^c$ and we have separated them by disjoint open sets.

Now, we want to convince ourselves that the space $\bN$ has the finite intersection property.

One direction of the finite intersection property is trivial, so let us assume that we are given an intersection $C = \bigcap_{\a\in\A} C_\a$ such that no finite intersection of closed sets $C_\a$ is empty. We need to show that the whole intersection is nonempty.

Since we are dealing with intersections only, without loss of generality we can assume that all the sets $C_\a$ are actually basic closed sets $\close A_\a$ with $A_\a\subseteq\N$. Now, the point is that the intersection of these sets is given by a filter: we have $C = \closed(\F)$ where the filter
\[ \F = \{ B\in\N : B \supseteq {A_{\a_1} \cap A_{\a_2} \cap \dots \cap A_{\a_n}} \text{ for some sets } A_{\a_i} \text{ with } \a_i\in\A\} \]
consists of all supersets of finite intersections of the sets $A_\a$. Since we assumed that none of these were empty, we have $\F\not\owns\emptyset$ and the collection $\F$ is indeed a filter. But we know that this filter corresponds to a \emph{nonempty} closed set $C=\closed(\F)$ by the previous Proposition \ref{prop-filter-closed}.
\end{proof}

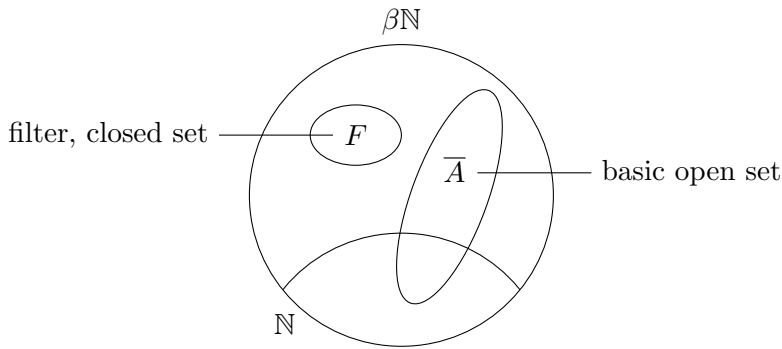
\begin{figure}[ht]
\begin{center}
\def\circleBetaN{(0,0) circle (2)}
\def\circleN{(0,-2.5) circle (2)}
\begin{tikzpicture}
	\draw \circleBetaN;
    \begin{scope}
        \clip \circleBetaN;
        \draw \circleN;
    \end{scope}
    % basic
	\node[anchor=south] at (90:2) {$\bN$};
    \node[anchor=north] at (-137:2.1) {$\mathbb{N}$};
    
    % basic open set
    \draw[rotate=-20,shift={(0.6,0.2)}] (0,0) ellipse (0.5 and 1.5);
	\node at (0.7,0.35) {$\close{A}$};
	\draw (1.0,0.3) -- +(1.5,0) node[right] {basic open set};
	% filter
	\draw (-0.6,0.8) ellipse (0.6 and 0.4) node {$F$};
	\draw (-0.9,0.8) -- +(-1.5,0) node[left] {filter, closed set};
\end{tikzpicture}
\end{center}
\caption{Illustration of the space of ultrafilters $\bN$. Basic open sets $\close A$ must contain natural numbers while closed sets $F$ need not.\label{pic-topology}}
\end{figure}

It is time to give a summary of the topology of $\bN$, illustrated by Figure \ref{pic-topology}.

%: FIGURE: topology of Stone-Cech compactification

The open sets, in particular the basic open sets of the form $\close A$, always contain ordinary numbers. They will be important because they allow us to transport properties of ultrafilters back to the natural numbers. This is another way to look at the permanence principle \ref{prop-permanence-principle}. Also note that sets containing only principal ultrafilters are automatically open.

The closed sets, on the other hand, do not necessarily contain natural numbers. Instead, they correspond to filters, i.e.\ arbitrary intersections of basic open sets. The Fréchet filter corresponds to the closed set of non-principal ultrafilters $\bN\setminus\N$.

The fact that filters correspond to \emph{nonempty} closed sets is closely related to the fact that the space $\bN$ is compact. Namely, the finite intersection property mimics the way filters behave: any finite intersection of sets in the filter is again in the filter and thus nonempty, which already implies that the intersection of all of its sets is nonempty as well. We cannot expect such an intersection to contain any ``old'' numbers from $\N$ though, that is too much to ask. For instance, the Fréchet filter corresponds to the intersection of sets $\close{\{n,n+1,\dots\}}$ which would have to contain an ``infinite natural number''\footnote{In this light, it is not surprising that ultrafilters feature prominently in the construction of Nonstandard Analysis, which deals with infinite and infinitesimal numbers. See also \textcite{Robinson:1996}.}.

%%%%%%%%%%%%%%%%%%%%%%%%%%%%%%%%%%%
\subsection{Limits along ultrafilters}
Many mathematicians encounter ultrafilters for the first time when they hear about Tychonoff's theorem, which belongs to point-set topology and states that arbitrary, even infinite products of compact topological spaces are again compact. This fundamental property of compactness is usually proven with ultrafilters, which are used to generalize \emph{sequences} and their \emph{convergence}.

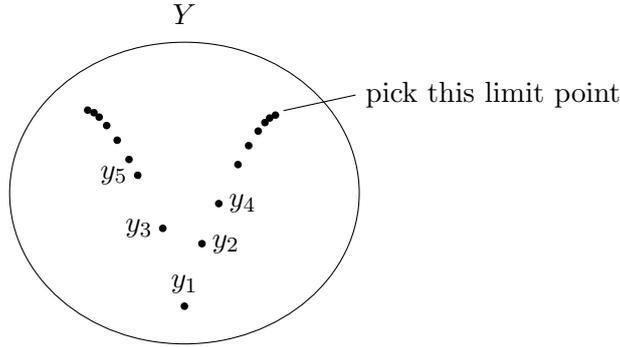
\begin{figure}[!ht]
\begin{center}
\begin{tikzpicture}
	\draw (0,0) ellipse (2.3 and 2);
	%\path (-130:2.2) node[below] {$Y$};
	\path (90:2.1) node[above] {$Y$};
	
	\def\mypoint{\fill circle (0.05)}
	\def\y#1{\mypoint node[above=1pt] {$y_{#1}$}}
	\def\yl#1{\mypoint node[left] {$y_{#1}$}}
	\def\yr#1{\mypoint node[right] {$y_{#1}$}}
	\def\yend{\mypoint node[above] (yend) {$y$}}
	
	\begin{scope}[decoration={markings
		,mark=at position 0.00 with \y1;
		,mark=at position 0.36 with \yl3;
		,mark=at position 0.62 with \yl5;
		,mark=at position 0.70 with \mypoint;
		,mark=at position 0.80 with \mypoint;
		,mark=at position 0.88 with \mypoint;
		,mark=at position 0.93 with \mypoint;
		,mark=at position 0.96 with \mypoint;
		,mark=at position 0.99 with \mypoint;
		}]
		\draw decorate{ (0,-1.5) .. controls (0,-0.8) and (130:1.6) .. (140:1.7) };
	\end{scope}
	\begin{scope}[decoration={markings
		,mark=at position 0.30 with \yr2;
		,mark=at position 0.5 with \yr4;
		,mark=at position 0.70 with \mypoint;
		,mark=at position 0.80 with \mypoint;
		,mark=at position 0.88 with \mypoint;
		,mark=at position 0.93 with \mypoint;
		,mark=at position 0.96 with \mypoint;
		,mark=at position 0.99 with \mypoint;
		}]
		\draw decorate{ (0,-1.5) .. controls (0,-1) and (50:1.5) .. (40:1.6) };
		\draw (40:1.7) -- (30:2.6) node[right] {pick this limit point};
	\end{scope}
\end{tikzpicture}
\end{center}
\caption{In a compact space $Y$, every sequence will have one or more limit points. An ultrafilter limit picks one of them.\label{pic-limit-point}}
\end{figure}

In that spirit, we are now going to define \emph{limits along ultrafilters}. The basic idea is the following: imagine a sequence $\sequence yn$ of elements in a compact topological space $Y$, for instance $Y=[0,1]$. In general, this sequence will not have a well-defined limit $\lim_{n\to\infty} y_n$, most likely because it oscillates wildly or does other strange things. But since the space $Y$ is \emph{compact}, we know at least that the sequence must have one or more \emph{limit points}, i.e.\ points $y\in Y$ that it comes close to infinitely often. Now, the idea is that ultrafilters allow us to attribute a well-defined limit to the sequence anyway, by arbitrarily but consistently choosing one of these limits points, as Figure \ref{pic-limit-point} visualizes.

%: FIGURE: limit points

In the following, we shall interpret sequences $\sequence yn$ as functions $f:\N\to Y, f(n) = y_n$ because this will give a more convenient notation.

%:	definition
\begin{definition}[Limit along an ultrafilter]
\label{definition-ultrafilter-limit}
Let $p\in \bN$ be an ultrafilter on the natural numbers and let $f:\N\to Y$ be a sequence of elements in a \emph{compact} topological space $Y$. Then, we define the \emph{limit of the sequence $f(n)$ along the ultrafilter}\index{limit along an ultrafilter|see{ultrafilter limit}}\index{ultrafilter limit} $p$ to be the unique point
\[  y =: \lim_{n\to p} f(n)  \in Y\]
which fulfills the condition
\[  \text{For every open neighborhood } U\owns y, \text{ we have } \inverse{f}(U)\in p  .\]
In other words, the sets of indices $\inverse f(U) = \{n\in\N : f(n) \in U\}$ of sequence elements that are mapped into the open neighborhood $U$ are required to be members of the ultrafilter $p$.
\end{definition}

For a principal ultrafilter $p$, this simply reduces to picking one member of the sequence, while for a non-principal ultrafilter $p$, this amounts to picking a limit point, because all the member sets of the ultrafilter are infinite.

The definition comes with a proof obligation, namely we have to show that such a limit $y$ actually exists and that it is unique. The compactness of the space $Y$ is crucial for this.

\begin{proof}[Limits along ultrafilters are well-defined]

\emph{Existence.}
\def\f{\tilde f}
The sequence $f:\N\to Y$ gives to rise to a \emph{map of ultrafilters} $p \mapsto \f(p)$ defined as
\[  \f(p) \owns A \defineiff p \owns \inverse{f}(A)  .\]
Since the inverse on sets $\inverse f$ preserves boolean operations, the collection $\f(p)$ is indeed an ultrafilter.

Now, choose the limit point $y$ to be some member of the intersection of all closed sets in the ultrafilter $\f(p)$,
\[  y \in \bigcap \{ A\subseteq Y : \f(p)\owns A, A \text{ closed }\} \]
This intersection is nonempty because $Y$, being compact, has the finite intersection property.

It remains to check that the point $y$ fulfills the defining condition for limits along ultrafilters, which the same as showing that the ultrafilter $\f(p)$ contains every open neighborhood of $y$. But if an open neighborhood $U\owns y$ were not a member of the ultrafilter $\f(p)$, then its complement $U^c$ would be both a \emph{closed set} and a member of $\f(p)$, which contradicts the construction of our point $y$.

\emph{Uniqueness.} This follows from the fact that $Y$ is Hausdorff. Namely, any two possible limit points $y$ and $y'$ can be separated by disjoint open neighborhoods $U_y$ and $U_{y'}$, but the ultrafilter $\f(p)$ can only contain one of these neighborhoods as their intersection is empty.
\end{proof}

The traditional approach for choosing limit points would be to pass to a subsequence. For instance, in functional analysis, it is not uncommon to consider a sequence in some compact space, like the unit ball of a Hilbert space in the weak topology, and then choose a convergent subsequence to get a limit point. However, this procedure can become very messy, especially when one has to repeat it. Limits along ultrafilters offer a clean alternative; in a sense, an ultrafilter has already picked all suitable subsequences ``in advance''.

The two approaches are not incompatible, though; it is always possible to turn subsequences into ultrafilters, as the following lemma demonstrates.

%We will make use of this in chapter \ref{chapter-measure} when we define counting measures on the space of ultrafilters $\bN$.

\begin{lemma}[Ultrafilters from subsequences]
\label{prop-subsequences}
Let $Y$ be a compact space and let $f:\N\to Y$ be a sequence. Furthermore, assume that a subsequence $(f(n_j))_{j=1}^\infty$ converges. Then, there exists an ultrafilter $p\in\bN$ that gives rise to the same limit for this particular sequence
\[  \lim_{n\to p} f(n) = \lim_{j\to\infty} f(n_j)   .\]
\end{lemma}
\begin{proof}
Any non-principal ultrafilter $p$ taken from the set of indices $\close A = \close{\{n_1,n_2,\dots\}}$ $\subseteq \bN$ will do the trick. After all, for any open neighborhood $U$ of the limit along the subsequence, we have
\[  \{n_k,n_{k+1},\dots\} \subseteq \inverse{f}(U)  \]
for some index $k\in\N$. But the set on the left-hand side is a member of $p$.
\end{proof}

Conversely, if we know that an ultrafilter belongs to a particular subsequence or set of indices, we may conclude that the ultrafilter limit only depends on the elements belonging to the subsequence.

\begin{lemma}[Restriction to a subset of indices]
\label{prop-limit-restriction}
Let $p$ be an ultrafilter contained in a particular basic open set $\close A\owns p$. Then, for any sequence $f:\N \to Y$ in a compact space $Y$, we have
\[  \lim_{n\to p} f(n) \in \closure(f(A))   .\]
In other words, the limit is determined only by the sequence elements with indices from the set $A\subseteq\N$.
\end{lemma}
\begin{proof}
If the limit did not lie in the closure $\closure(f(A))$, then it would be contained in an open set $U\subseteq Y$ that does not intersect $f(A)$. But then, we would have $\inverse f(U) \subseteq A^c \not\in p$, a contradiction to the definition of the ultrafilter limit.
\end{proof}

%:	rules for calculation

We still have to clarify what is meant by picking a limit point ``consistently''. It simply means that the limit along ultrafilter commutes with continuous functions, as we would expect of any well-behaved notion of ``limit''. Also, all the usual rules for calculating with limits retain their validity. Here just a sample.

\begin{proposition}[The usual limit rules apply to ultrafilter limits]
\index{ultrafilter limits! rules for calculation}
\label{ultrafilter-limit-addition}
\label{ultrafilter-limit-rules}
Let $X$ be a compact topological space and let $f:\N\to X$ be a sequence. Then, the following statements hold:
\begin{enumerate}
\item Let $h:X\to Y$ be a continuous function into another compact space $Y$. Then, we may interchange limits 
\[  \lim_{n\to p} h(f(n)) = h (\lim_{n\to p} f(n))  .\]
\item Let $g:\N\to X'$ be a sequence in another compact space $X'$. Then, we may interchange limits and pairs
\[  \lim_{n\to p} (f(n),g(n)) = \left(\lim_{n\to p} f(n), \lim_{n\to p} g(n)\right)  .\]
\item Assume that $X=[0,1]$ is a compact interval of real numbers and let $g:\N\to X$ be another sequence. Then, we can interchange limits and summation
\[  \lim_{n\to p} (f(n)+g(n)) = \lim_{n\to p} f(n) + \lim_{n\to p} g(n)  .\]
\end{enumerate}
\end{proposition}

Proving these rules directly would be rather tedious, not to mention that there are many more of them. It is much more efficient to recast the notion of ultrafilter limit in terms of a universal property, giving rise to the so-called \emph{\StoneCech}. The rules then follow from a simple uniqueness argument.

%:	Stone-Cech compactification

\label{stone-cech}
\begin{definition}[\StoneCech]
\label{definition-stone-cech}
The \define{\StoneCech\ $\beta X$} of a topological space $X$ is a compact topological space together with a continuous map $\iota: X\to\beta X$ subject to the following universal property: any map $f : X \to Y$ into another compact space $Y$ factors uniquely through a map $\beta f:\beta X \to Y$, as the following diagram indicates
\[\begin{tikzpicture}[scale=1.5,auto]
	\node (bX) at (0,1) {$\beta X$};
	\node (X) at (0,0) {$X$};
	\node (Y) at (2,0) {$Y$};
	
	\draw [->] (X) edge node {$\iota$} (bX);
	\draw [->] (X) edge node[below] {$f$} (Y);
	\draw [->] (bX) edge node[pos=0.2,right=14pt] {$\exists!\,\beta f$} (Y);
\end{tikzpicture}\]
\end{definition}

It is a standard exercise to show that the \StoneCech\ is actually unique up to homeomorphisms.

As our choice of notation suggests, the space of ultrafilters $\bN$ is indeed the \StoneCech\ of the natural numbers $\N$. This is true for every discrete space $X$, but not for general topological spaces. In the latter case, there exist different ultrafilters that cannot be distinguished by any continuous function $f$. An extreme example would be a compact space $X$, because then we already have $\beta X = X$.

\begin{proof}[$\bN$ is the \StoneCech\ of $\N$]
\def\f{\beta f}

\emph{Existence.}
Of course, the map $\iota$ is the standard embedding and the map $\f$ is given by the limit along ultrafilters:
\[  \f(p) :=  \lim_{n\to p} f(n).  \]
We have to show that it is continuous. To that end, let $U$ be any open neighborhood of a point $\f(p)$ in the target space $Y$.

First, let us construct a smaller open neighborhood $V \owns \f(p)$ with the property that its closure is contained in the original one, $\closure V \subseteq U$. Doing that is a standard exercise in point-set topology: use the compactness of the space $Y$ to construct two disjoint open sets $V$ and $W$, i.e.\ $V \cap W = \emptyset$, that separate the point $\f(p) \in V$ from the closed set $U^c \subseteq W$. The set $V$ will have the desired property since the set $W^c$ is closed and hence $V \subseteq \closure V \subseteq W^c \subseteq U$.

Now, consider the basic open set $\close{\inverse{f}(V)}$. By definition of the ultrafilter limit, this is a basic open neighborhood of the ultrafilter $p$. But Lemma \ref{prop-limit-restriction} about restrictions to subsets of indices implies that it is mapped into the open set $U$,
\[  \f(\close{\inverse{f}(V)}) \subseteq \closure(f(\inverse{f}(V))) = \closure(V) \subseteq U .\]
This proves that the map $\f$ is continuous.

\emph{Uniqueness.} Let $g : \bN\to Y$ be another continuous map that makes the diagram commute. We have to show that
\[  g(p) = \lim_{n\to p} f(n)  .\]

Since the map $g$ is continuous, we know that the preimage of any open neighborhood $U$ of the point $g(p)$ contains a basic open set, $\inverse{g}(U) \supseteq \close A\owns p$. But our maps $f$ and $g$ agree on the natural numbers, so $A\subseteq \inverse{f}(U)\in p$. In other words, the point $g(p)$ fulfills the very definition of the limit in question.
\end{proof}

Now, this universal property easily gives the rules for calculating with ultrafilter limits.

\begin{proof}[of Proposition \ref{ultrafilter-limit-rules}, Rules for calculating with ultrafilter limits]
\begin{enumerate}
\item Consider the diagram
\[\begin{tikzpicture}[scale=1.5,auto]
	\node (bN) at (0,1) {$\bN$};
	\node (N) at (0,0) {$\N$};
	\node (X) at (2,0) {$X$};
	\node (Y) at (4,0) {$Y$};
	
	\draw [->] (N) edge node {} (bN);
	\draw [->] (N) edge node[below] {$f$} (X);
	\draw [->] (X) edge node[below] {$h$} (Y);
	\draw [->] (bN) edge node[pos=0.7,left=14pt] {$\beta f$} (X);
	\draw [->] (bN) edge node[pos=0.2,right=28pt] {$\beta (h\circ f)$} (Y);
\end{tikzpicture}\]
Uniqueness of the rightmost diagonal arrow implies that $\beta(h\circ f) = h \circ \beta f$, which is just another way of writing the equation
\[ \lim_{n\to p} h(f(n)) = h (\lim_{n\to p} f(n)) .\]
\item The product of continuous functions $\beta f \times \beta g$ makes the diagram
\[\begin{tikzpicture}[scale=1.5,auto]
	\node (bN) at (0,1) {$\bN$};
	\node (N) at (0,0) {$\N$};
	\node (X) at (2,0) {$X\times X'$};
	
	\draw [->] (N) edge (bN);
	\draw [->] (N) edge node[below] {$f\times g$} (X);
	\draw [->] (bN) edge node[pos=0.2,right=14pt] {$\beta f\times \beta g$} (X);
\end{tikzpicture}\]
commute. But by uniqueness, this means that $\beta(f \times g) = \beta f \times \beta g$, which is again merely another way of writing the equation to be proven.
\item Decompose the sum of maps $f+g = h\circ (f \times g)$ as a pair of maps followed by the continuous function $h(x,x') = x+x'$ and apply the two previous statements.
\end{enumerate}
\end{proof}

The \StoneCech\ also allows us to identify the Banach space of bounded sequences $\boundedSequences$ with a Banach space continuous functions, namely $\Continuous(\bN)$. We will make use of this important fact when constructing \emph{measures} on the space of ultrafilters $\bN$ in Section $\ref{section-measure-construction}$.

\begin{theorem}[Bounded sequences as continuous functions on $\bN$]
\label{prop-bounded-sequences}
There is a canonical isomorphism of Banach spaces
\[  \boundedSequences \xrightarrow{\sim} \Continuous(\bN) \]
given by the ultrafilter limit
\[ \sequence fn \mapsto \left(p \mapsto \lim_{n\to p} f_n\right).\]
\end{theorem}
\begin{proof}
Any sequence $f\in\boundedSequences$ bounded by $|f_n| \leq C$ can be viewed as a function $f:\N\to[-C,C]$ into a compact interval. It lifts to a continuous function $\beta f: \bN \to [-C,C]$ which is given by the ultrafilter limit.

Conversely, since the \StoneCech\ $\bN$ is compact, any continuous function $g:\bN\to\R$ must be bounded, $|g(p)| \leq C$. The restriction to natural numbers $g_n := g|_{\N}(n)$ gives rise to a bounded sequence.

Thanks to the rules for calculating with limits, this identification preserves vector space operations. It is also easy to check that the Banach norms are, in fact, preserved. Hence, we have an isomorphism of Banach spaces.
\end{proof}

%\to do{Maybe example: Banach limits for summing divergent series?}
% No, they are boring. Only bounded series. The really divergent series have unbounded terms.

%%%%%%%%%%%%%%%%%%%%%%%%%%%%%%%%%%%%%%%%%%%%%%%%%%%%%%%%%%%%%%%%%%%%%%
\section{Ramsey theory and ultrafilters}
\label{chapter-ramsey}
%%%%%%%%%%%%%%%%%%%%%%%%%%%%%%%%%%%
\subsection{What is Ramsey theory?}
\index{Ramsey theory}
Our interest in ultrafilters is actually motivated by an interest in \emph{Ramsey theory}, a branch of mathematics which we will now describe. It is named after the British mathematician Frank P. Ramsey (1903-1930).

The following theorem is a prototypical example of a statement from Ramsey theory.
\begin{theorem}[van der Waerden]
\label{prop-van-der-Waerden}
\index{van der Waerden's theorem}
Consider the set of natural numbers $\N$ and imagine that we paint each number with one of $r$ different colors. In other words, consider a partition
\[ \N = \partition C \]
where $C_i$ is the subset of natural numbers painted with the $i$-th color. Then, for any given length $k$, there exists at least one arithmetic progression
\[ a,a+b,\dots,a+(k-1)b \qquad \text{with }b>0 \]
that is \emph{monochromatic}, i.e.\ whose members all have the \emph{same color}. Put differently, for any length $k$, there is one color $i$ such that the set $C_i$ contains an arithmetic progression of length $k$.
\end{theorem}

We already know that the set of natural numbers $\N$ contains arithmetic progressions aplenty, but what van der Waerden's theorem says is that no matter how we partition this set into finitely many pieces, one of them will \emph{also} contain an arithmetic progression. This surprising phenomenon merits a proper name.

\begin{definition}[Partition regularity]
Let $\mathcal G$ be a nonempty collection of sets that we deem interesting, for instance
\[  \mathcal G = \{ \{a,a+b,\dots,a+(k-1)b\} : a,b\in\N, b > 0\}  \]
the collection of arithmetic progressions of length $k$. Such a collection is called \define{partition regular} if, for every partition
\[  \N = \partition{C}, \quad C_i\subseteq\N  ,\]
at least one of the parts $C_i$ still contains an interesting set $G\in\mathcal G$, i.e.\ $G \subseteq C_i$.

To make language more convenient, we will often talk about the \emph{property} ``contains an interesting set $G\in\mathcal G$'' being partition regular, instead of the \emph{collection} $\mathcal G$ being partition regular.
\end{definition}

In other words, van der Waerden's theorem asserts that the property ``contains an arithmetic progression of length $k$'' is partition regular.

In general, Ramsey theory is the study of partition regularity: take a mathematical object, like $\N$, or some graph, and divide it into finitely many parts. Then, provided that the object in question is sufficiently large, one of the parts will always contain an interesting substructure, in our case an arithmetic progression.

Usually, the reason for containing interesting substructures is that one of the parts is also ``large'' in a suitable sense, so large that it cannot possibly miss the interesting substructure. As we will see in Section \ref{section-ramsey-ultrafilter}, being a member of an ultrafilter is such a suitable notion of largeness. What other notions of largeness can guarantee that a set of natural numbers contains arithmetic progressions?

One example is the notion of a \emph{syndetic set}.

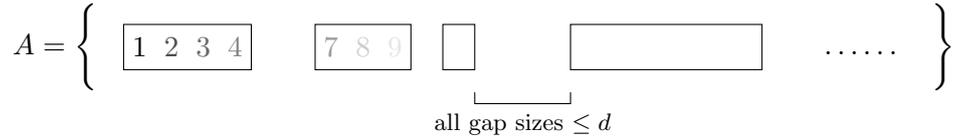
\begin{figure}[ht]
\begin{center}
\begin{tikzpicture}[yscale=0.5,xscale=0.7,scale=0.6]
	\def\block#1#2{\draw (#1,-1) rectangle +(#2,2)}
	\def\number#1[#2]{\node[#2] at (#1+0.5,0) {$#1$}}

	%\draw[help lines] (0,-2) grid (10,2);
	\node[left] at (0.5,0) {$A = \Bigg\{$};
	\block{1}{4};
	\number{1}[black];
	\number{2}[black!70];
	\number{3}[black!60];
	\number{4}[black!50];
	\block{7}{3};
	\number{7}[black!40];
	\number{8}[black!20];
	\number{9}[black!10];
	\block{11}{1};
	\block{15}{6};
	
	\draw (12,-2) -- +(0,-0.5) -- node[midway,below,style={font=\footnotesize}] {all gap sizes $\leq d$} (15,-2.5) -- (15,-2);
	
	\node at (25,0) {$\dots\dots \quad \Bigg\}$};
\end{tikzpicture}
\end{center}
\caption{Illustration of a syndetic set. The rectangular blocks indicate natural numbers that are contained in the set, while empty space indicates numbers that are absent from the set. Being syndetic means that no gap may exceed a fixed size $d$.\label{figure-syndetic}}
\end{figure}

\begin{definition}[Syndetic]
A set $A$ of natural numbers is called \define{syndetic} if it has \emph{bounded gaps}. This means that there exists a size $d\in\N$ such that every interval $[M,M+d]\subseteq\N, M\in\N$ has nonempty intersection with the set $A$. Put differently, this implies that adjacent elements from the set $A$ are no more than a distance $d$ away from each other, as illustrated by Figure \ref{figure-syndetic}.
\end{definition}

%: FIGURE: syndetic set

As you can see, syndetic sets contain at least a fraction $1/d$ of all natural numbers, so it is certainly appropriate to call them ``large''. There is an equivalent definition of syndeticity which, together with van der Waerden's theorem, shows that syndetic sets also contain arithmetic progressions.

\begin{proposition}[Syndetic, alternate definition]
A set $A\subseteq\N$ is syndetic if and only if the set of natural numbers $\N$ can be written as a union of finitely many shifts of $A$,
\[  \N = \bigcup_{i=1}^k (A - n_i) , \quad\text{for some }n_1,\dots,n_k\in\N.\]
\end{proposition}
\begin{definition}[Shift]
\label{definition-shift}
Let $A\subseteq\N$ be a set of natural numbers. For a number $n\in\N$, the \emph{shift}\index{shift of a set} $A-n$ is defined as the set $A-n := \{ k\in\N : k+n\in A\}$, or more suggestively
\[ k\in A-n \defineiff (k+n) \in A.\]
In other words, subtract $n$ from the numbers in the set $A$ and keep the nonnegative results.
\end{definition}
\begin{proof}[Syndetic, alternate definition]
If $A\subseteq\N$ is a syndetic set with maximal distance $d\in\N$, then the union $\bigcup_{i=1}^d (A-i)$ is clearly equal to $\N$. On the other hand, for a given collection of shifts $(A-n_1),\dots,(A-n_k)$, we can simply choose $d=\max\{n_1,\dots n_k\}$.
\end{proof}
\begin{corollary}[Syndetic sets contain arithmetic progressions]
Let $A\subseteq\N$ be a syndetic set. Then, for every length $k\in\N$, the set $A$ contains an arithmetic progression of length $k$.
\end{corollary}
\begin{proof}
A finite union of shifts $A-n_i$ covers the set $\N$. Of course, van der Waerden's theorem applies not only to paritions but also to covers, so we obtain an arithmetic progression of length $k$ in one of the shifts $C_i:=A-n_i$. But we can simply add $n_i$ to each member of the progression and obtain an arithmetic progression in the set $A$.
\end{proof}

Another notion of largeness for a set, perhaps the most natural one, is to contain a positive fraction of all natural numbers, to have a \emph{positive density}. The density of a set $A$ is determined by the ratios $|A\cap [1,N]|/N$ as $N\to\infty$. For example, the set $A=2\N$ of even numbers has density $1/2$ and is large in this sense.

Unfortunately, these ratios do not converge in general, which leads to a variety of different definitions of density. We will be concerned with the following variant, which no longer requires the interval to begin at the number $1$:

\begin{definition}[Upper Banach density]
\label{definition-upper-Banach-density}
The \define{upper Banach density} of a set $A\subseteq\N$ of natural numbers is defined to be the limes superior
\[ \d(A) := \limsup_{N\to\infty}\sup_{M\in\N} \frac{|A\cap [M+1,M+N]|}{N} \]
taken over shifted intervals $[M+1,M+N]$ whose length goes to infinity.
\end{definition}

Is there any reason to expect that a set with positive upper Banach density contains arithmetic progressions? Probably not, but a deep theorem first proven by Szeméredi, previously conjectured by Erdős and Turán, asserts that this is nonetheless true.

\begin{theorem}[Szemerédi]
\index{Szemeredi's theorem@Szemerédi's theorem}
\label{theorem-szemeredi}
Every set of natural numbers with positive \density{} contains arithmetic progressions of every length.
\end{theorem}

This is a formidable generalization of van der Waerden's theorem. To see that, note that the \density{} is subadditive,
\[ \d(A\cup B) \leq \d(A) + \d(B), \quad\text{for any }A,B\subseteq\N.\]
Thus, any partition of the natural numbers gives rise to the inequality
\[ 1 = \d(\N) = \d(\partition C) \leq \d(C_1) + \d(C_2) \dots + \d(C_r)\]
which implies that at least one of the sets $C_i$ must have positive \density{} and hence contain arithmetic progressions of every length.

We will make no attempt to prove this theorem here; for that, you are referred to the original proof by \textcite{Szemeredi:1975p1591} or the classical proof relying on ergodic theory by \textcite{Furstenberg:1980p1531}. Instead, our goal is to extend the interpretation of van der Waerden's theorem in terms of ultrafilters, to be presented in Section \ref{section-ramsey-ultrafilter}, to Szemerédi's theorem, which will be done in Section \ref{section-szemeredi-ultrafilter}.

%%%%%%%%%%%%%%%%%%%%%%%%%%%%%%%%%%%
\subsection{Partition regular properties correspond to ultrafilters}
\label{section-ramsey-ultrafilter}
We now want to explain the connection between partition regularity and ultrafilters.

The key insight is foreshadowed by the following trivial observation: the property ``contains the number 7'' is partition regular. After all, if we partition the set $\N$ into finitely many parts, one of them must contain the number $7$. Unlike, say, an arithmetic progression, the set $\{7\}$ cannot be broken apart because it consists of just a single point.

But we have seen that the natural numbers are not the only points of $\N$: the ultrafilters also behave like ``points''. Could it be that van der Waerden's theorem holds because there exists a ``magical point'' $p$ such that any set containing this point automatically contains an arithmetic progression?

This is indeed the case! Namely, we now show that every partition regular property is \emph{equivalent} to the existence of ultrafilters whose member sets all have the property in question. One could quip: ``Ramsey theory is just the search for funny ultrafilters.''

\begin{theorem}[Equivalence of partition regularity and ultrafilters]
Let $\mathcal G$ be a collection of sets that we deem interesting. Then, the following are equivalent:
\begin{itemize}
\item $\mathcal G$ is partition regular.
\item There exists an ultrafilter $p$ such that all member sets $A\in p$ contain at least one interesting set $G\in\mathcal G$, i.e.\ $G\subseteq A$.
\end{itemize}
\end{theorem}
\begin{proof}
``$\Longleftarrow$'': The ultrafilter $p$ is the ``magical point'' we were talking about. Namely, consider a partition $\N = \partition C$ and take the \emph{closure} of this partition in the space of ultrafilters $\bN$. Since closure commutes with boolean operations, we have
\[  \bN = \close{\partition C} = \partition{\close C} .\]
Clearly, one of the parts $\close C_i$ must contain the ultrafilter $p$, which is just another way of saying that $C_i\in p$. By definition of the ultrafilter $p$, this implies that the set $C_i$ contains an interesting set $G\in\mathcal G$, i.e.\ $G\subseteq C_i$. Hence, the collection $\mathcal G$ is partition regular.

%: FIGURE: partition of \beta\N
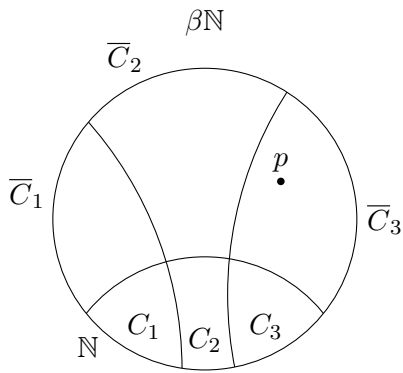
\begin{figure}[ht]
\begin{center}
\def\circleBetaN{(0,0) circle (2)}
\def\circleN{(0,-2.5) circle (2)}
\begin{tikzpicture}
	\draw \circleBetaN;
    \begin{scope}
        \clip \circleBetaN;
        \draw \circleN;
        
        % division
        \draw (-5.3,-2) circle (5);
        \draw (5.3,-1) circle (5);

		% point p        
        \fill (1,0.5) circle (0.05) node[above] {$p$};
    \end{scope}
    \node[above] at (90:2.3) {$\bN$};
    \node[below] at (-137:2.1) {$\mathbb{N}$};
    
    \node[left] at (170:2) {$\close C_1$};
    \node[above] at (120:2.1) {$\close C_2$};
    \node[right] at (0:2) {$\close C_3$};
    
    \node[above] at (-115:1.9) {$C_1$};
    \node[above] at (-90:1.9) {$C_2$};
    \node[above] at (-65:1.9) {$C_3$};
\end{tikzpicture}
\end{center}
\caption{A partition $\N = \partition C$ extends seamlessly to the space of ultrafilters $\bN$. Clearly, one of the parts must contain the ultrafilter $p$.}
\end{figure}
%% Proof

``$\Longrightarrow$'': Given a collection $\mathcal G$ which is partition regular, we have to construct an ultrafilter $p$ whose members sets all contain interesting subsets. Of course, we are going to use the ultrafilter construction lemma (\ref{prop-ultrafilter-construction}) for that.

Our predicate will be the notion of partition regularity applied to any subset of the natural numbers $\N$:
\begin{align*}
 \phi(A) := \text{``}&\text{For every partition $A=\partition C$,} \\
  & \text{one of the parts contains an interesting set $G\subseteq C_i, G\in\mathcal G$.''}
\end{align*}
By assumption, $\phi(\N)$ holds true, so our starting point will be the trivial filter $\mathcal F := \{\N\}$. It is also easy to see that $\phi(A)$ fulfills the second requirement.

For the third requirement, we argue by contradiction. Let $A=A_1\uplus A_2$ be a disjoint union and assume that $\phi(A)$ holds true while both $\phi(A_1)$ and $\phi(A_2)$ are false. In other words, there exist partitions
\[ A_1 = \partition B \quad\text{and}\quad A_2 = \partition C \]
such that none of the parts contain an interesting set. But clearly, this gives rise to a joint partition
\[ A = \partition B \uplus \partition C \]
such that none of the parts contain an interesting set, which would imply that $\phi(A)$ is false, a contradiction.

Hence, our predicate fulfills all the requirements and there exists an ultrafilter $p$ with the property that $\phi(A)$ for all member sets $A\in p$. But clearly, the property $\phi(A)$ implies that the set $A$ contains an interesting subset, just consider the trivial partition $A=A$. 
\end{proof}

In other words, van der Waerden's theorem can be thought of as a theorem about the existence of an ultrafilter whose members all contain arithmetic progressions, instead of as a theorem about partitions of the set $\N$.

Note, however, that this ultrafilter is by no means unique: a single partition regular property may have many different ultrafilters associated to it. Let us introduce a notation for these ultrafilters in the case of arithmetic progressions.

\begin{definition}[Progression-rich ultrafilters]
\label{definition-aprich}
An ultrafilter $p\in\bN$ is called \emph{\aprich}\index{ultrafilter!progression-rich} if every member set $A\in p$ contains an arithmetic progression of length $k\in\N$. The set of \aprich\ ultrafilters for a particular length $k$ is denoted with
\[  AP_k := \left\{p\in\bN : \begin{array}{l}\text{all } A\in p \text{ contain at least}\\[-5pt]\text{one arithmetic progression of length } k
\end{array}\right\} .\]
The set of \aprich\ ultrafilters that feature every length is denoted with
\[  AP_\infty := \bigcap_{k=1}^\infty AP_k  .\]
\end{definition}

Thus, if the closure $\close A$ of a set of natural numbers contains a \aprich\ ultrafilter $p\in AP_k$, then the set $A$ is guaranteed to contain an arithmetic progression of length $k$. The \emph{converse} is not true, however: there may be sets of natural numbers which contain many arithmetic progressions of length $k$ but are not members of some ultrafilter $p\in AP_k$. That is because containing an ultrafilter is a stronger property, it is also stable under partitions $A = A_1 \uplus A_2$. The situation is much better for the set $AP_\infty$, though; in the next section, we will argue that a set contains arithmetic progressions of every length exactly when it is a member of some ultrafilter $p\in AP_\infty$.

To conclude this section, we give a somewhat explicit description of the sets $AP_k$ that will become important when we prove van der Waerden's theorem in Section \ref{proof-van-der-Waerden}. Also, it will show that these sets are closed, which we will use when we link them to Szemerédi's theorem in Section \ref{section-szemeredi-ultrafilter}. Note that, in conjunction with the finite intersection property, this also implies that the set $AP_\infty$ is nonempty as long as all the sets $AP_k$ are nonempty.

\begin{proposition}[Sets of \aprich\ ultrafilters as preimages]
\label{prop-apk}
\label{prop-aprich}
Let $\bN^k$ denote the $k$-fold product of the topological space $\bN$ with itself and denote the diagonal map by
\[  \diagonalMap_k : \bN \to \bN^k, x \mapsto (x,x,\dots,x)  .\]
Furthermore, consider the set of ``arithmetic'' tuples
\[  AT_k := \{(a,a+b,\dots,a+(k-1)b) : a,b\in\N,b > 0\} \subseteq\N^k\subset\bN^k.\]
Then, the preimage of its closure in the space $\bN^k$ is precisely the set of \aprich\ ultrafilters
\[ AP_k = \inverse{\diagonalMap_k}(\closure_{\bN^k}(AT_k)) \subseteq \bN. \]
In particular, the set of \aprich\ ultrafilters is closed because it is the preimage of a closed set under a continuous function.
\end{proposition}
\begin{proof}
\def\p{\tilde p}\def\B{\tilde B}
Let $p\in\bN$ be an ultrafilter and $\p=\diagonalMap_k(p)=\diagonal p$ its image.

By definition, a point $\p$ is contained in the closure $\closure_{\bN^k}(AT_k)$ if and only if every one of its rectangular open neighborhoods $\B = \fold {\close A}\times{k}$ intersects the set $AT_k$. But since $\p$ is in the diagonal, we can restrict ourselves to quadratic neighborhoods of the form $\B = \close A\times \close A\times\dots\times\close A$ where $\close A=\fold {\close A}\cap{k}$ is a basic open neighborhood of the ultrafilter $p$. Such a set $\B$ intersects the set of arithmetic tuples $AT_k$ exactly if the set $\close A$ contains an arithmetic progression of length $k$.

To summarize, the point $\p$ is contained in the closure $\closure_{\bN^k}(AT_k)$ if and only if all basic open neighborhoods $\close A$ of the ultrafilter $p$ contain an arithmetic progression of length $k$.
\end{proof}

%%%%%%%%%%%%%%%%%%%%%%%%%%%%%%%%%%%
\subsection{Finitary statements}
We have introduced partition regularity in terms of partitions of the infinite set $\N$, but Ramsey theory can also have a different, more ``finitary'' flavor. Let us now also take a brief look at the finitary point of view and the interesting compactness argument that translates between both views. In particular, it will allow us to characterize the ultrafilters $p\in AP_\infty$ specifically in terms of their member sets.

Originally, van der Waerden proved the following statement.

\begin{proposition}[van der Waerden, finitary version]
For each count $r$ of colors and each length $k$, there exists a natural number $W(k,r)$, called the \define{van der Waerden number}, such that for every size $N\geq W(k,r)$, any coloring of the interval
\[  [1,N] = \partition C \]
will contain a monochromatic arithmetic progression of length $k$.
\end{proposition}

Here, the emphasis is on the size $W(k,r)$ of the interval: if you choose it large enough, then one of its parts will always contain an interesting substructure. It is an interesting but difficult question to give good bounds for the numbers $W(k,r)$, see the overview by \textcite{Green:2007} for references to recent results.

It is clear that this version of van der Waerden's theorem implies our previous version \ref{prop-van-der-Waerden} that concern partitions of the infinite set $\N$. But the converse is also true, as the following argument shows.

\begin{proof}[van der Waerden, finitary version follows from infinitary version]

Assuming that every coloring of the set $\N$ with $r$ colors will contain monochromatic arithmetic progressions of length $k$, we have to show that the same is already true of a finite interval $[1,W(k,r)]$ with large enough size $W(k,r)\in\N$.

We argue by contradiction: Let us assume that there exists a sequence of growing intervals $[1,N^n]$ and partitions $[1,N^n]=C_1^n \uplus C_2^n \uplus \dots \uplus C_r^n$ that do \emph{not} contain a monochromatic arithmetic progression of length $k$. We have to use this sequence to construct a partition of the set of natural numbers $\N = \partition C$ that does not contain monochromatic arithmetic progressions either.

Each partition from the sequence can be represented by a function
\[ c^n : [1,N^n] \to \{1,\dots,r\},\quad c^n(k) = i \quad\text{if and only if } k \in C^n_i \]
that labels each number $k$ with its color $i$. To construct the global partition, choose an arbitrary non-principal ultrafilter $p$ and consider the function given by the ultrafilter limit
\[ c : \N\to\{1,\dots,r\},\quad c(k) := \lim_{n \to p} c^n(k) .\]
The limit exists because the target space $\{1,\dots,r\}$ is compact and the values $c^n(k)$ are well-defined as soon as the interval size $N^n$ becomes large enough. Another way to interpret this construction is that the \emph{space of partitions} $\{1,\dots,r\}^\N$ is compact and that our sequence of partitions has a limit point.

It remains to be shown that the limit point $c$ does not contain a monochromatic arithmetic progression either. But this condition is essentially a collection of local properties. For instance, consider the arithmetic progression $\{a,a+d,a+2d\}$ for some numbers $a,d\in\N,d>0$. Its elements have different colors if and only if
\[ |c(a) - c(a+d)| + |c(a+d) - c(a+2d)| + |c(a+2d) - c(a)| \geq 1  .\]
But this inequality was true for the functions $c^n$ and is preserved in the ultrafilter limit. Hence, the partition represented by the function $c$ does not contain any monochromatic arithmetic progression.
\end{proof}

As an application, we can now prove that any set which contains arithmetic progressions of every length must be a member of some ultrafilter $p\in AP_\infty$. Of course, this result is conditional on the infinitary version of van der Waerden's theorem, which we will prove only much later in Section \ref{section-ideals}.

\begin{proposition}[Ultrafilters describing long arithmetic progressions]
A set $A$ of natural numbers contains arithmetic progressions of every length if and only if its closure $\close A$ contains an ultrafilter $p\in AP_\infty$.
\end{proposition}
\begin{proof}

\Backward: This direction is immediate from the definition of the set $AP_\infty$ (Definition \ref{definition-aprich}).

\Forward: Consider the predicate
\[ \phi(B) = \text{``The set $B$ contains arithmetic progressions of every length''}.\]
We will show that it fulfills the conditions from the ultrafilter construction lemma (\ref{prop-ultrafilter-construction}). To that end, choose $\mathcal F$ to be the principal filter on the set $A$, $\mathcal F = \{F \subseteq \N : A \subseteq F\}$. It is clear that the predicate fulfills the first and second requirement. To prove the third requirement, we have to show that for any set $B$ with $\phi(B)$ and any partition $B = B_1 \uplus B_2$, we have $\phi(B_1)$ or $\phi(B_2)$ as well.

Let $k\in\N$ be a length. By assumption, the set $B$ contains an arithmetic progression $P$ of length $K:=W(k,2)$. Of course, the partition $B = B_1 \uplus B_2$ induces a partition $P=P_1\uplus P_2$ of the arithmetic progression. Since the collection of arithmetic progressions is invariant under affine maps, we can apply the finitary version of van der Waerden's theorem to the \emph{progression} $P$ instead of the \emph{interval} $[1,W(k,2)]$ and obtain that one of the sets $P_1\subseteq B_1$ or $P_2\subseteq B_2$ must contain an arithmetic progression of length $k$.

Hence, we have shown that for any length $k$, at least one of the sets $B_1$ or $B_2$ must contain an arithmetic progressions this length. By the pigeonhole principle, this means that one of the sets must contain arithmetic progressions for infinitely many lengths $k$, which already implies that it contains them for all lengths $k$. In other words, we have $\phi(B_1)$ or $\phi(B_2)$ as desired.

Now, the ultrafilter construction lemma yields an ultrafilter $p$ with $\phi(B)$ for every $B \in p$, i.e.\ $p\in AP_\infty$. But our choice of the filter $\mathcal F$ also gives $A \in p$ as desired.
\end{proof}

%%%%%%%%%%%%%%%%%%%%%%%%%%%%%%%%%%%%%%%%%%%%%%%%%%%%%%%%%%%%%%%%%%%%%%
\section{Algebra in \titlebN}
\label{chapter-algebra}
%%%%%%%%%%%%%%%%%%%%%%%%%%%%%%%%%%%
\subsection{Addition of ultrafilters, limit version}
As we have seen in Section \ref{section-ramsey-ultrafilter}, van der Waerden's theorem is equivalent to the existence of an ultrafilter whose member sets all contain arithmetic progressions. We now want to prove van der Waerden's theorem by constructing this ultrafilter directly, instead of deducing its existence from an elementary, but intricate combinatorial proof.

Our trusted ultrafilter construction lemma (\ref{prop-ultrafilter-construction}) is no longer useful for this purpose. Instead, we turn towards a fascinating algebraic structure on the space of ultrafilters, namely \emph{addition} of ultrafilters.

\begin{definition}[Addition of ultrafilters, limit version]
The \emph{sum}\index{ultrafilter sum}\index{ultrafilter addition|see{ultrafilter sum}} of two ultrafilters $p,q\in\bN$ is defined as the double limit
\[  q + p := \left(\lim_{m\to q}\ \lim_{n\to p}\ (m + n)\right)  \in\bN ,\]
of the ordinary sum $m+n$ of natural numbers $m,n\in\N$.
This limit is well-defined because the topological space $\bN$ is compact.
\end{definition}

Unfortunately, the addition of ultrafilters no longer has the same nice properties as the addition of ordinary numbers that we are used to. For instance, it is not commutative, and it is not even a continuous operation on the space $\bN$. The main reason for that is that the two ultrafilter limits cannot be interchanged,
\[  \lim_{m\to q}\ \lim_{n\to p} \neq \lim_{n\to p}\ \lim_{m\to q} .\]
This is not very surprising: even ordinary limits cannot be interchanged at times, and ``forcing'' a limit to exist in the ultrafilter sense does not ameliorate this tendency.

Hence, when performing calculations, we may only rely on the following rules.

\begin{proposition}[Properties of the addition of ultrafilters]
Addition of ultrafilters is
\begin{enumerate}
\item \define{left-continuous} (``continuous in the left argument''), which means that the map
\[  \rho_p : \bN \to \bN,\quad \rho_p(q) = q+p \]
%\qquad\text{(``add the element $p$ on the right'')}
is continuous for every ultrafilter $p\in\bN$. However, addition is not continuous in both arguments. Instead, only for ordinary numbers $m\in\N$, we can be assured that the map
\[  \lambda_m : \bN \to \bN,\quad \lambda_m(p) = m+p \]
is continuous as well.
\item not commutative, but the natural numbers $\N$ lie in the \emph{center}. In other words, for every natural number $m\in\N$, we have
\[ p+m = m+p\quad\text{ for all } p\in\bN .\]
\item \emph{associative}, i.e.\
\[  r + (q + p) = (r + q) + p  \quad\text{ for all } p,q,r\in\bN  .\]
\end{enumerate}
\end{proposition}
\begin{proofx}
\begin{enumerate}
\item Proving continuity is equivalent to checking that the functions commute with ultrafilter limits. By definition of the sum, for any natural number $m\in\N$, we have
\[  \lim_{n\to p} \lambda_m(n) = \lim_{n\to p}(m + n) =: m+p = m + \lim_{n\to p} n = \lambda_m(p).\]
This proves the second claim. Using this, we obtain the first claim
\[  \rho_p(q) = q + p := \lim_{m\to q}\lim_{n\to p}(m + n) = \lim_{m\to q} (m+p) = \lim_{m\to q} \rho_p(m).\]
\item Continuity allows us to use the commutativity of the natural numbers $\N$. Namely, thanks to the previous property, for any natural number $m\in\N$, we have
\[ p+m = \lim_{n\to p} \lim_{m\to m}(n+m) = \lim_{n\to p} (n+m) = \lim_{n\to p} (m+n) = m + \lim_{n\to p} n = m + p .\]
\item Use the previous properties to slide natural numbers past the limits:
\begin{align*}
  r + (q + p) &= \lim_{k\to r} (k + \lim_{m\to q} \lim_{n\to p} (m+n)) = \lim_{k\to r} \lim_{m\to q} \lim_{n\to p} (k + m+n ) \\
  &= \lim_{k\to r} (\lim_{m\to q} \lim_{n\to p} (k + m) + n) = (r + q) + p
.\end{align*}\qedEquationHack{0.98cm}
\end{enumerate}
\end{proofx}

We did not prove that commutativity really fails for non-principal ultrafilters, but we will make no use of this fact. Consult \textcite{Hindman:1998p531} for a proper proof.

In subsequent sections, we also want to consider slightly extended versions of ultrafilter addition, like component-wise addition of tuples from the space $\bN^k$, or restrictions to closed subsets thereof. Let us unify them under the term \emph{left topological semigroup}.

\begin{definition}[Left topological semigroup]
Let $S$ be a topological space and $(+) : S\times S\to S$ be a binary operation. The pair $(S,+)$ is called a \define{left topological semigroup} if the operation $(+)$ is
\begin{properties}
\item associative, and
\item left-continuous (``continuous in the left argument''), which means that the function
\[ \rho_p : S \to S,\quad \rho_p(q) = q + p \]
is continuous in the argument $q$ for any element $p\in S$.
\end{properties}
\end{definition}

In other words, $(\bN,+)$ is a left topological semigroup. Other examples of interest to us are the space of tuples $(\bN^k,+)$ with component-wise addition and left topological semigroups $(M,+)$ coming from subsets $M\subseteq \bN^k$ with $M+M\subseteq M$.

Since our semigroup operations are not commutative, it might not be a good idea to denote them with the symbol ``$+$''; but the author thinks that the notation is justified because all of our example are extensions of the addition of natural numbers.

Another interesting example of a left topological semigroup is the set of non-principal ultrafilters $\bN\setminus\N$. Remember that this set is closed, hence compact. But why is the sum of two non-principal ultrafilters again non-principal? The following lemma will help us to give a proof by clarifying the relationship between closure and ultrafilter addition.

Note that we use the standard notation for sumsets: $A+B = \{ a+b : a\in A, b \in B\}$.

\begin{lemma}[Closure and addition in left topological semigroups]
\label{prop-addition-closure}
Let $S$ be a left topological semigroup and let $A,B\subseteq S$ be two sets. Then, the following statements are true.
\begin{enumerate}
\item $(\closure A) + B \subseteq \closure(A+B)$
\item If the semigroup $S$ is compact, then we have
\[ (\closure A) + b = \closure (A+b) \]
for every element $b\in B$.
\item If the set $A$ commutes with every element of the semigroup $S$, then we have
\[ (\closure A) + (\closure B) \subseteq \closure(A + B) .\]
For instance, the set $A$ might be a set of natural numbers in the semigroup $S=\bN$.
\end{enumerate}
\end{lemma}
\begin{proofx}
\begin{enumerate}
\item Addition to the right $\rho_b(a) = a + b$ is continuous. Hence, for every element $b\in B$, the preimage $\inverse{\rho_b}(\closure (A + b))$ is closed. Since it contains the set $A$, it also contains the closure $\closure A$, which means $(\closure A) + b \subseteq \closure(A + b)$. Taking the union over all elements $b\in B$ gives the result.

\item Using the previous statement, we only have to show that the set $(\closure A) + b$ is closed. But the set $\closure A$ is compact and addition to the right maps compact sets to compact sets.

\item Remembering the first statement and using commutativity, we get
\begin{align*}
(\closure A) + (\closure B) &\subseteq \closure(A + \closure B) \subseteq \closure(\closure B + A)
\\ &\subseteq \closure(\closure(B + A)) = \closure(A + B) .
\end{align*}
\end{enumerate}\qedEquationHack{0.98cm}\vskip\parskip
\end{proofx}

\begin{proposition}[Sum of non-principal ultrafilters is non-principal]
Writing $\corona\N = \bN\setminus\N$ for the set of non-principal ultrafilters, we have
\[ \corona\N + \corona\N = \corona\N .\]
\end{proposition}
\begin{proof}
\def\Segment#1{\N_{\geq #1}}
Let $\Segment{k} := \{k,k+1,…\}$ be the set of natural numbers greater or equal than the number $k$. We have already seen that the set of non-principal ultrafilters is the intersection of their closures, $\corona\N = \bigcap_{k=1}^\infty \close{\Segment k}$.

The main ingredient to the proof is the observation that for any numbers $n,k\in\N$, we have
\[ n + \Segment k \subseteq \Segment k .\]
Hence, thanks to the previous lemma, we have
\[ n + \corona\N = n + \bigcap_{k=1}^\infty \close{\Segment k} = \bigcap_{k=1}^\infty \close{n + \Segment k} \subseteq \bigcap_{k=1}^\infty \close{\Segment k} = \corona\N ,\]
which even implies $\bN + \corona\N \subseteq \corona\N$.
\end{proof}

%%%%%%%%%%%%%%%%%%%%%%%%%%%%%%%%%%%
\subsection{Addition of ultrafilters, set membership version}
While the definition of ultrafilter sums in terms of limits is easy to calculate with, it does not give a good description of the ultrafilter $p+q$ in terms of its member sets $A\in p+q$. To that end, we want give an alternative description of ultrafilters sums, which is often used as a definition. But first, let us introduce some handy notation.

\begin{definition}[Ultrafilter shift]
\label{definition-ultrafilter-shift}
Let $p\in\bN$ be an ultrafilter and $A\subseteq \N$ be a set. We say that the \define{ultrafilter shift} $A-p$ is the set of natural numbers defined by 
\[  n \in A-p \defineiff A-n\in p .\]
\end{definition}

This notation for ultrafilter shifts is taken from \textcite{Beiglbock:2009p581}. Remember that we can regard ordinary numbers $m\in\N$ as principal ultrafilters via the equivalence $B \in m \iff m \in B$. Then, the ultrafilter shift $A-m$ will simply be the usual shift, see Definition \ref{definition-shift}.

Keep in mind that the shift $A-p$ is a set of \emph{natural numbers}, unlike the sumset $A+p$, which would be a set of \emph{ultrafilters}. To avoid ambiguity, we will never use the minus sign to denote sets of ultrafilters. Fortunately, making a mistake with this rule is not a serious issue; it is an easy exercise to deduce ``$\close{A - p} = \close A - p$'' from the next proposition.

In any case, with this notation, we can now describe the member sets of the sum of two ultrafilters.
\begin{proposition}[Addition of ultrafilters, set membership version]
\label{prop-addition-membership}
The member sets of a \emph{sum}\index{ultrafilter sum} $q+p$ of two ultrafilters $p,q\in\bN$ are given by the equivalence
\[  A \in q + p \iff A - p \in q  .\]
\end{proposition}

The notation is very suggestive: simply ``subtract'' $p$ from the right on both sides of the member relation $\in$. We will collect a few more rules like this in a moment and obtain a small but efficient calculus for manipulating ultrafilter membership. Then, we can prove many theorems by almost mechanically shifting ultrafilters back and forth the membership sign $\in$.

Note that these rules work best when using equivalences to reason about set membership and avoiding the curly braces. For instance, we could equally have defined the ultrafilter shift with the notation $A-p := \{ n \in \N : A - n \in p\}$, but in the author's opinion, this form is neither easy to decode mentally nor is it conductive to formal manipulation. Having said that, curly braces can be handy at times, for instance for denoting preimages as in the following proof; but it seems that they are best manipulated with just the trivial rule $\{n : n\in A\} = A$.

\begin{proof}[Addition of ultrafilters, equivalence of the two definitions]
Consider a basic open set $\close A\subseteq\bN$. Unraveling the definition of the ultrafilter limit twice, we obtain the following expression with many curly braces:
\begin{align*}
A \in q+p &\iff \lim_{m\to q}\lim_{n\to p}\ (m + n) \in \close A \iff \{ m\in\N: \lim_{n\to p}\ (m + n) \in \close A\} \in q
\\&\iff \{n\in\N : \{ m\in\N : m+n \in\close A\} \in p\} \in q
.\end{align*}
Since $m+n$ is always a natural number, we can simplify the inner set to
\[ \{ m\in\N : m+n \in\close A\} = \{ m\in\N : m+n \in A\} = \{ m\in\N : m \in A-n\} = A-n .\]
Applying the definition of the ultrafilter shift gives
\[ \{n\in\N : A-n \in p\} = \{n\in\N : n \in A-p \} = A-p.\]
and hence
\[ A \in q+p \iff A-p \in q \]
as desired.
\end{proof}

We now list and prove the syntactic rules for easy manipulation of ultrafilter sums and their member sets. They are to be used in conjunction with the other rules, like the definition of closure \ref{definition-closure}, the permanence principle \ref{prop-permanence-principle} and the rules for boolean algebra \ref{prop-boolean-algebra}.

\label{prop-syntactic-rules}
\begin{proposition}[Syntactic rules for ultrafilter sums and shifts]
\index{ultrafilter sum!syntactic rules}
Let $A\subseteq\N$ be a set of natural numbers, $p,q,r\in\bN$ be ultrafilters and $n,m\in\N$ be natural numbers. Then, the following laws hold.
\begin{align}
A\in q+p &\iff A-p \in q  	 \tag{Addition}\\
A\in p+n &\iff A\in n+p \tag{Commutes with $\N$}\\
(A-p)-q &= A - (q+p)\footnotemark{}	\tag{Associativity} \\
(A\cap B) - p &= (A-p) \cap (B-p) \tag{Shift of intersection}\\
(A\cup B) - p &= (A-p) \cup (B-p) \tag{Shift of union}
\end{align}
\footnotetext{The reversal of summands is due to non-commutativity and similar in appearance to the rule $\inverse{(gh)} = \inverse h\inverse g$ for groups.}
\end{proposition}
\begin{proofx}
\begin{enumerate}
\item We have already shown this.
\item We already know this, but there is also a direct proof using the definition of the ultrafilter shift:
\[A \in n+p \iff A-p \in n \iff n \in A-p \iff A-n \in p \iff A \in p+n .\]
\item For every $n\in\N$, associativity tells us that
\begin{align*}
n \in (A-p)-q &\iff A-p \in n+q \iff A \in (n+q)+p
\\&\iff A \in n+(q+p) \iff n \in A-(q+p) 
.\end{align*}
\item For all $n\in\N$, Proposition \ref{prop-ultrafilters-points} on ultrafilters behaving like points allows us to infer
\begin{align*}
n \in (A\cap B) - p &\iff (A\cap B) \in n+p \iff A\in n+p \text{ and } B\in n+p
\\&\iff (A-p)\in n \text{ and } (B-p)\in n
\\&\iff (A-p)\cap(B-p) \owns n
.\end{align*}
\item Similar to the previous rule.\hfill\qedsymbol
\end{enumerate}
\end{proofx}

One last remark on syntax and notation shall conclude this section. Namely, when defining ultrafilter addition, we could equally well have chosen to take the right limit last instead of the left limit, with the consequence that the role of the left and right argument would have been switched and that we would now be talking about ``right topological semigroups''. This is the convention that some authors, like \textcite{Bergelson:2003p505} adopt, while other authors, like \textcite{Hindman:1998p531}, agree with our convention. It may seem that for reasons of symmetry, one notation is as good as the other, and readers will have to suffer eternal confusion unless mathematicians agree to arbitrarily prefer one direction over the other. However, it appears that the symmetry of the situation is, in fact, broken by the minus sign ``$-$''. The author thinks that the notation $A-p \in q \iff A \in q + p$ should be preferred over its mirror image ``$(-p) + A \in q \iff A \in p + q$'' and this leads to the convention adopted here.\footnote{Note that this does not settle the question whether addition should be called ``left continuous'' (``continuous in its left argument'') or ``right continuous'' (``adding something to the right is continuous''), but at least we would agree on the order of arguments.}

%%%%%%%%%%%%%%%%%%%%%%%%%%%%%%%%%%%
\subsection{Idempotent ultrafilters}
As we have seen in the previous section, the addition $+$ of natural numbers can be extended to the space of ultrafilters $\bN$, but many of the known rules, like commutativity, are lost. In general, analyzing the algebraic structure of the left topological semigroup $(\bN,+)$ is a difficult task, in particular because the axiom of choice gives us only a cursory grasp of its elements. And if you ask the wrong question, like how many elements of a certain kind does it contain, you will soon run into foundational issues that are independent of the ZFC axioms of set theory.

Hence, we can only ask and answer rather simple questions; one of the simplest questions being whether there exist \define{idempotent} ultrafilters, i.e.\ elements $p\in\bN$ with the property that $p+p=p$. Fortunately, these simple questions have nontrivial consequences when translated into statements about partition regularity, which is the very beauty of the ultrafilter approach. In particular, the existence of idempotent ultrafilters implies an interesting theorem by Hindman (Corollary \ref{theorem-hindman}) about sets of finite sums and they also enable us to prove van der Waerden's theorem for the special case of $k=3$ (Theorem \ref{prop-van-der-Waerden-k-3}).

It is not difficult to prove that \emph{finite} semigroups always contain idempotent elements, so it is not surprising that this extends to the case of \emph{compact} left topological semigroups like $\bN$ or $\bN^k$, as we will now show.

\begin{theorem}[Ellis. Existence of idempotent ultrafilters]
\index{Ellis' theorem}
\label{prop-Ellis}
Let $(M,+)$ be a compact left topological semigroup; for instance a closed set of tuples $M\subseteq \bN^k$ with $M+M\subseteq M$. Then, the semigroup $M$ always contains an \define{idempotent} element, i.e.\ an element $p\in M$ with the property $p+p=p$.
\end{theorem}
\begin{proof}
Consider the collection of the nonempty closed (and compact) subsets $N\subseteq M$ which also fulfill $N+N \subseteq N$. We want to apply Zorn's lemma to obtain a set $N$ which is \emph{minimal} with respect to inclusion. Clearly, the intersection $\bigcap_\alpha N_\alpha$ over any chain $N_1 \supseteq N_2 \dots$ of sets from the collection is again in the collection. In particular, since the space $M$ is compact and has the finite intersection property, this intersection is \emph{nonempty}. We claim that our minimal set $N$ contains the desired element $p$.

First, consider the set $N+p$ for some arbitrary element $p\in N$. By continuity, this set is compact. Associativity implies that $(N+p) + (N+p) \subseteq N+p \subseteq N$. But since the set $N$ was chosen to be minimal, we have
\[ N = N+p .\]

In particular, there exists an element $q\in N$ such that $q+p = p$. Now, consider the set of all such elements
\[ L := \{ q \in N : q + p = p \} .\]
Associativity implies that $L+L \subseteq L$ because
\[ q,q'\in L \implies (q' + q) + p = q' + (q + p) = q' + p = p \implies (q'+q) \in L.\]
Furthermore, the set $L$ is closed, because it is the preimage of the closed set $\{p\}$ under the continuous function $q \mapsto q+p$.

Appealing to the minimality of our set $N$, we conclude that $L=N$. This means $p\in L$, which implies $p+p=p$ as desired.
\end{proof}

In particular, applying Ellis' theorem to the left topological semigroup $\bN\setminus\N$ yields the existence of non-principal idempotent ultrafilters.

Translated into Ramsey theory, such non-principal idempotent ultrafilters give rise to Hindman's theorem, which says that the property of containing a \emph{set of finite sums} is partition regular. Let us define what that means.

\begin{definition}[Set of finite sums]
Let $\sequence nk$ be an increasing sequence of natural numbers $n_k\in\N$. Its \define{set of finite sums} $FS(\sequence nk)$ consists of all the finite sums of the sequence elements
\[\begin{split}
FS(\sequence nk) &:= \left\{\sum_{k\in F} n_k : F \text{ a finite subset of } \N\right\} \\
  &= \left\{n_{k_1} + n_{k_2} + \dots + n_{k_l} : k_i\in\N, k_1 < k_2 < \dots < k_l\right\}
\end{split}\]
\end{definition}

Sets containing sets of finite sums are called \emph{IP-sets}; the abbreviation ``IP'' stands for ``IdemPotent''.

\begin{definition}[IP-set]
An \define{IP-set} is a set $A\subseteq\N$ that contains a set of finite sums, $FS(\sequence nk) \subseteq A$.
\end{definition}

\begin{theorem}[Ultrafilters characterizing IP-sets]
A set $A$ of natural numbers is an IP-set if and only if its closure $\close A$ in the space $\bN$ contains a non-principal idempotent ultrafilter.
\end{theorem}
\begin{proof}
\Backward :
Let $p\in\close A$ be an idempotent ultrafilter. Since $p=p+p$, we have
\[ A \in p \iff A \in p + p \iff A-p \in p .\]
In particular, we also have
\[ A \cap (A-p) \in p .\]
By the permanence principle \ref{prop-permanence-principle}, there exists a natural number $n_1\in\N$ contained in this intersection, i.e.\ $n_1\in A$ and $n_1 \in A-p$. The latter relation is equivalent to $A-n_1\in p$, which implies
\[ A\cap(A-n_1) \in p .\]
Repeating this procedure with the set $A\cap(A-n_1)$ and using the permanence principle to choose a natural number $n_2\in A$ with $n_2 > n_1$, we obtain
\[ (A\cap(A-n_1))\cap[(A-n_2)\cap((A-n_1)-n_2)] \in p .\]
This can be continued indefinitely and we obtain a sequence of natural numbers $(n_1 < n_2 < …)$ such that
\[ \{n_1, n_2, n_1+n_2, …\} = FS(\sequence nk) \subseteq A \]
as desired.

It may be instructive to sketch a different, though weaker argument that uses idempotence more directly. Namely, consider the ultrafilter limit of a finite tuple
\[ \lim_{n_1\to p} \lim_{n_2\to p} \lim_{n_3\to p} (n_1,n_2,n_3,n_1+n_2,n_1+n_2+n_3) = (p,p,p,p+p,p+p+p) .\]
Due to idempotence, this is equal to the diagonal tuple $(p,p,p,p,p)$. Similarly to the characterization of \aprich\ ultrafilter in Proposition \ref{prop-aprich}, this implies that every basic open neighborhood of the ultrafilter $p$ must contain all elements of a finite tuple.

\Forward :
Assume that the set $A$ contains a set of finite sums $FS(\sequence nk) \subseteq A$ and consider the following subsets
\[  S_l := FS((n_k)_{k=l}^{\infty}) \subseteq FS(\sequence nk) \subseteq A  \]
which only contain sums of elements $n_k$ with indices starting from $k \geq l$. By the finite intersection property, the intersection of their closures
\[  S := \bigcap_{l=1}^\infty \close S_l  \quad\subseteq \close A \cap (\bN\setminus\N)\]
is nonempty. We claim that this is a semigroup, $S+S=S$. Then, by Ellis' theorem, it must contain the desired non-principal idempotent ultrafilter.

The key to showing that the set $S$ is a semigroup is the following observation:
\[ \text{For every sum } t\in S_l, \text{ there exists an index } m \geq l \text{ such that } t+S_m \subseteq S_l .\]
After all, the sum $t$ has the form $t = n_{k_1} + n_{k_2} + … + n_{k_\tau}$ with a greatest index $k_\tau$, and adding to that any sum that starts with the index $m := k_\tau + 1$ will again give a proper sum from the set $S_l$.

Since $t$ is an ordinary number and commutes with all ultrafilters, we can conclude that
\[ t + S_m \subseteq S_l \implies t + \close S_m \subseteq \close S_l \implies t + S \subseteq \close S_l .\]
Choosing the sum $t\in S_l$ arbitrarily and taking the intersection over all indices $l$, we obtain
\[ \forall l.\ S_l + S \subseteq \close S_l \implies \forall l.\ \close S_l + S \subseteq \close S_l \implies S + S \subseteq S \]
as desired.
\end{proof}

Translating this characterization of IP-sets into Ramsey theory, we obtain the following theorem as a corollary.
\begin{corollary}[Hindman]
\index{Hindman's theorem}
\label{theorem-hindman}
Let $\N = \partition{C}$ be a partition of the natural numbers. Then, at least one of the parts $C_i$ is an IP-set, i.e.\ contains a set of the form $FS(\sequence nk)$.
\end{corollary}

Amusingly, idempotent ultrafilters also enable us to prove a special case of van der Waerden's theorem, namely for arithmetic progressions of length $3$.

\begin{theorem}[van der Waerden, $k=3$]
\index{van der Waerden's theorem!special@special case $k=3$}
\label{prop-van-der-Waerden-k-3}
Let $p$ be any idempotent ultrafilter. Then, every open neighborhood of the ultrafilter $2p+p$ contains an arithmetic progression of length $k=3$.
\end{theorem}

Here, the multiplication of an ultrafilter by a natural number is defined as
\[ 2p := \lim_{n\to p} 2n .\]
In particular, this operations is continuous. Note that $2p \neq p+p$ in general!

\begin{proofx}[van der Waerden, $k=3$]
Remembering the characterization of \aprich\ ultrafilters from Proposition \ref{prop-apk}, we have to show that the diagonal triple $(2p+p,2p+p,2p+p) \in \bN^3$ is contained in the closure of the set of arithmetic triples
\[  AT_3 := \{(a,a+b,a+2b) : a,b\in\N,b > 0\} \subseteq\bN^3 .\]
We will achieve this by expressing it as the ultrafilter limit of a sequence in the closure $\closure AT_3$. Namely, we claim that
\[  (2p+p,2p+p,2p+p) \?= \lim_{a\to 2p} \lim_{b\to p} \lim_{c\to p} (a+c,a+b+c,a+2b+c) \in \closure AT_3.\]
Evaluating the right-hand side with the rules \ref{ultrafilter-limit-rules} for calculating with ultrafilter limits gives
\begin{align*}
 \text{r.h.s}
  &= \lim_{a\to 2p} \lim_{b\to p} ( \lim_{c\to p} (a+c), \lim_{c\to p} (a+b+c), \lim_{c\to p} (a+2b+c)) \\
  &= \lim_{a\to 2p} \lim_{b\to p} (a+p,a+b+p,a+2b+p) \\
  &= \lim_{a\to 2p} (a+p,a+p+p,a+2p+p) \\
  &= (2p+p,2p+p+p,2p+2p+p)
.\end{align*}
We obtain the desired left-hand side by using that the ultrafilter $p$ is idempotent and simplifying $p+p=p$ and $2p+2p=2(p+p)=p$. Of course, the latter equation requires justification:
\begin{align*}
2p+2p
  &= \lim_{m\to p} 2m + \lim_{n\to p} 2n = \lim_{m\to p}\lim_{n\to p} (2m + 2n) \\
  &= \lim_{m\to p}\lim_{n\to p} 2(m+n) = 2\lim_{m\to p}\lim_{n\to p} (m+n) = 2(p+p)
.\end{align*}
\qedEquationHack{0.98cm}
\end{proofx}

%%%%%%%%%%%%%%%%%%%%%%%%%%%%%%%%%%%
\subsection{Ideals and minimal ideals}
\label{section-ideals}
We now want to prove van der Waerden's theorem in full generality by exhibiting a \aprich\ ultrafilter. To that end, we will study the \emph{ideals} of the left topological semigroup $\bN$ and see that ultrafilters contained in the \emph{minimal left ideals} are always \aprich.

\begin{definition}[Notions of ideals]
\index{ideal!left, right, two-sided}
\index{ideal!minimal left ideal}
\index{ideal!minimal two-sided ideal@minimal two-sided ideal $\kbN$}
Consider a left topological semigroup $(S,+)$, for example the space of ultrafilters $S=\bN$. A subset $I\subseteq S$ is called
\begin{enumerate}
\item a \emph{left ideal} if $p+I \subseteq I$ for every element $p\in S$.
\item a \emph{right ideal} if $I+p \subseteq I$ for every element $p\in S$.
\item a \emph{minimal left ideal} if it is a left ideal and does not contain a proper subset that is also a left ideal.
\item an \emph{ideal}, or \emph{two-sided ideal}, if it is both a left and a right ideal.
\item a \emph{minimal ideal} if it is a two-sided ideal and does not contain a proper subset that is also a two-sided ideal.
\end{enumerate}
\end{definition}

To keep things concrete, we will formulate the subsequent propositions about ideals for the semigroup $S = \bN$ only, even though they are valid for any compact left topological semigroup.

First, we present a collection of useful trivia about minimal left ideals. Being minimal, they have a fairly simple structure. Note, however, that a semigroup may contain many different minimal left ideals.

\begin{proposition}[Trivia about minimal left ideals]
\index{ideal!minimal left ideal}
\begin{enumerate}
\item A left ideal $L\subseteq \bN$ is minimal if and only if it is generated by each of its elements, i.e.\ if it has the form
\[  L = \bN + p \text{ for each } p\in L .\]
\item Minimal left ideals are compact.
\item Every left ideal $L\subseteq\bN$ contains at least one minimal left ideal.
\item Let $L\subseteq \bN$ be a minimal left ideal and consider any sum $q+p \in L$ of an element $p\in L$ and some ultrafilter $q\in\bN$. Then, there always exists another ultrafilter $t\in\bN$ such that
\[  t + q + p = p  .\]
\end{enumerate}
\end{proposition}
\begin{proofx}
\begin{enumerate}
\item For each element $p\in L$, the set $\bN + p$ is clearly a left ideal contained in $L$. If $L$ is minimal, then it must already be the whole of $L$.

Conversely, consider any left ideal $L'$ contained in $L$. Again, the set $\bN + p \subseteq L'$ for some $p \in L'$ is also a left ideal. But the assumption says that it is already equal to $L$, hence $L'=L$ and the left ideal $L$ must be minimal.

\item Sets of the form $\bN+p$ are closed.

\item Apply Zorn's lemma to the collection of closed left ideals contained in the set $L$. This collection is nonempty because for any $p\in L$, the left ideal $\bN+p \subseteq L$ is closed.

\item If the left ideal $L$ is minimal, then $p \in L = \bN + (q+p)$. \hfill\qedsymbol\vskip\parskip
\end{enumerate}
\end{proofx}

For minimal two-sided ideals, the situation is much nicer: There is only one minimal two-sided ideal.

\begin{proposition}[The minimal two-sided ideal $\kbN$]
\label{definition-minimal-ideal-kappa}
\index{ideal!minimal two-sided ideal@minimal two-sided ideal $\kbN$}
There exists exactly one minimal two-sided ideal, called $\kbN$. In other words, this ideal is contained in every other two-sided ideal.

Moreover, the set $\kbN$ is the union of all minimal left ideals,
\[ \kbN = \Union_{p\ \in\text{ minimal left ideal}} (\bN+p).\]
\end{proposition}
\begin{proof}
Let $\kbN$ be the union of all minimal left ideals as the formula indicates. Clearly, this set is a left ideal.

First, we have to show that it is also a right ideal. Let $\bN+p \subseteq\kbN$ be a minimal left ideal. We want to prove that for any ultrafilter $q\in\bN$, the left ideal $L := (\bN+p)+q$ is also minimal. Let $t+p+q\in L$ be an arbitrary element. By minimality of the left ideal $\bN+p$, there exists an ultrafilter $s\in\bN$ with $s+t+p = p$. But associativity implies $s+(t+p+q) = (p+q)$ and the left ideal $L$ is generated by $t+(p+q)$. In other words, the left ideal $L$ is generated by each of its elements, so it must be minimal.

Finally, we have to show that every two-sided ideal $I\subseteq\bN$ contains the set $\kbN$, which is the same as showing that the set $I$ contains every minimal left ideal $L$. Since $I$ is a right ideal and $L$ a left ideal, we have $\emptyset \neq I+L \subseteq I\cap L$, so the intersection of the two is definitely nonempty. But the intersection $I\cap L \subseteq L$ is a left ideal and by minimality of $L$, it must already be the whole of $L$. Hence, the left ideal $L=I\cap L$ is contained in the ideal $I$.
\end{proof}

\label{kappa-is-aprich}
\label{proof-van-der-Waerden}
We are now ready to prove van der Waerden's theorem by showing that ultrafilters from the minimal ideal $\kbN$ are \aprich.

\begin{theorem}[Ultrafilters in the minimal ideal are \aprich]
\label{prop-kbN-aprich}
\index{van der Waerden's theorem}
\index{ultrafilter!progression-rich}
Let $p\in\kbN$ be an element of the minimal ideal. Then, the ultrafilter $p$ is \aprich, $p\in AP_\infty$.
\end{theorem}
\begin{proof}
\def\p{\tilde p}\def\q{\tilde q}\def\t{\tilde t}
Remembering the characterization of \aprich\ ultrafilters from Proposition \ref{prop-apk}, we have to show that for every natural number $k$, the diagonal tuple $\p := \diagonalMap_k(p) = \diagonal p \in\bN^k$ is contained in the closure $M := \closure_{\bN^k} AT_k$ of the set of arithmetic $k$-tuples
\[ AT_k = \{(a,a+b,\dots,a+(k-1)b) : a,b\in\N,b > 0\} \subseteq\N^k\subset\bN^k  .\]

First, note that the set $AT_k$ is a semigroup, $AT_k+AT_k \subseteq AT_k$, and that we can also add diagonal tuples of numbers without leaving it, $AT_k + \diagonalMap_k(\N) \subseteq AT_k$. Since these sets lie in the center of the semigroup $\bN^k$, Lemma \ref{prop-addition-closure} tells us that taking closures preserves these relations. In other words, we obtain
\[ M + M \subseteq M \text{ as well as } \p + M \subseteq M \text{ and } M + \p \subseteq M .\]
Put differently, the set $M + \p$ is a closed semigroup containing a two-sided ideal $M$.

By Ellis' theorem \ref{prop-Ellis} on the existence of idempotents, the compact left topological semigroup $M+\p$ contains an element $\q+\p = (q_1+p,q_2+p,…,q_k+p)$ with $\q
\in M$ and the property
\[  (\q+\p) + (\q+\p) = \q + \p  .\]
Since the ultrafilter $p\in\bN$ was assumed to be contained in a minimal left ideal of the semigroup $\bN$, there exists a tuple of ultrafilters $\t = (\vector tk)\in\bN^k$ such that
\[ \t + \q + \p = \p  .\]
Now, if the tuple $\t$ were an element of $M$, then we would be done already, because this would imply $\t+\q\in M$ and $\p \in M + \p \subseteq M$. Unfortunately, this is not the case, but having chosen the element $\q+\p$ to be idempotent allows us to perform a somewhat similar argument.

Namely, we have
\[  \p + (\q + \p) = (\t + \q + \p) + (\q + \p) = (\t + \q + \p) = \p  ,\]
which implies $\p \in \p + M + \p \subseteq M$. This completes the proof.
\end{proof}

%We conclude this section with a remark on the ``Hales-Jewett-Theorem'', which is a generalization of van der Waerden's theorem from natural numbers $\N$ to free semigroups with more than one generator. To learn more about it, consult \textcite{Hindman:1998p531}. Here, we just want to mention that its proof is virtually identical to the one given above, except that one has to replace $\N$ with a free semigroup. \to do{Remove this rather unilluminating remark?}

%%%%%%%%%%%%%%%%%%%%%%%%%%%%%%%%%%%
\subsection{Syndetic sets and ultrafilters}
As we have seen, all ultrafilters in the minimal ideal $\kbN$ are \aprich, i.e.\ their member sets contain arithmetic progressions of every length. But the member sets are even more special than this: We will now show that they are exactly the \emph{piecewise syndetic} sets.

\begin{definition}[thick, piecewise syndetic]
A set of natural numbers is called
\begin{itemize}
\item \define{thick} if it contains arbitrarily large intervals.
\item \define{piecewise syndetic} if it is equal to the intersection of a thick set and a syndetic set.
\end{itemize}
\end{definition}

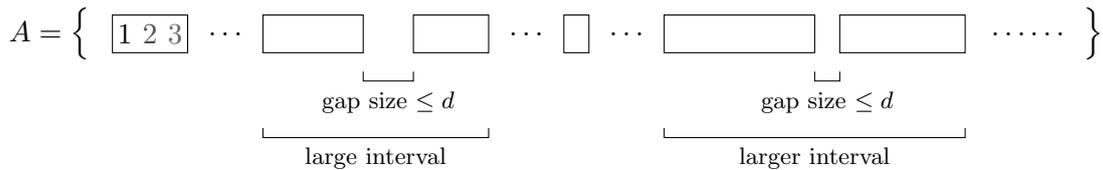
\begin{figure}[ht]
\begin{center}
\begin{tikzpicture}[yscale=0.5,xscale=0.66,scale=0.5]
	\def\block#1#2{\draw (#1,-1) rectangle +(#2,2)}
	\def\number#1[#2]{\node[#2] at (#1+0.5,0) {$#1$}}

	%\draw[help lines] (0,-2) grid (10,2);
	\node[left] at (0.5,0) {$A = \Big\{$};
	\block{1}{3};
	\number{1}[black];
	\number{2}[black!70];
	\number{3}[black!60];
	
	\node at (5.6,0) {$…$};
	
	% interval
	\block{7}{4};
	\block{13}{3};
	
	\def\gapsize#1#2{
	\draw (#1,-2) -- +(0,-0.5) -- node[midway,below,style={font=\footnotesize}] {gap size $\leq d$} (#2,-2.5) -- +(0,0.5)}
	\gapsize{11}{13};

	\def\intervals#1#2#3{
	\draw (#1,-5) -- +(0,-0.5) -- node[midway,below,style={font=\footnotesize}] {#3} (#2,-5.5) -- +(0,0.5)}
	\intervals{7}{16}{large interval};
	
	% small block
	\node at (17.6,0) {$…$};
	\block{19}{1};
	
	% interval
	\node at (21.6,0) {$…$};
	\block{23}{6};
	\block{30}{5};
	\gapsize{29}{30};
	\intervals{23}{35}{larger interval};
	
	\path (37.6,0) node {$\dots\dots$} +(2.5,0) node {$\Big\}$};
\end{tikzpicture}
\end{center}
\caption{Piecewise syndetic sets always consist of eventually growing intervals that have bounded gaps. The space between these intervals may be arbitrary, however.\label{figure-piecwise-syndetic}}
\end{figure}

Put differently, a piecewise syndetic set may contain arbitrarily large gaps, but there are always larger and larger intervals on which the set does have bounded gaps as illustrated by Figure \ref{figure-piecwise-syndetic}.

%: FIGURE: piecewise syndetic set

\begin{theorem}[Ultrafilters characterizing piecewise syndetic sets]
A set $A$ of natural numbers is \emph{piecewise syndetic} if and only if its closure $\close A$ contains an ultrafilter $p\in\kbN$ in the minimal ideal.
\end{theorem}

Translated into Ramsey theory, this implies that the property of being piecewise syndetic is partition regular.

\begin{figure}[h]
\begin{center}
\def\circleBetaN{(0,0) circle (2)}
\def\circleN{(0,-2.5) circle (2)}
\begin{tikzpicture}
	\draw \circleBetaN;
    \begin{scope}
        \clip \circleBetaN;
        \draw \circleN;
        
        % kbN
        \draw (0,0.7) ellipse (1 and 0.5);
        \begin{scope}
        	\clip (0,0.7) ellipse (0.9 and 0.45);
			\foreach \y in {-0.4,-0.2,0,0.2,0.4} {
				\draw (-1,0.7+\y) -- (1,0.7+\y);
			}
        \end{scope}
        
        % thick set
        \draw (0,-1) ellipse (0.8 and 2);
        
        % syndetic set
        \draw (5.3,-0.3) circle (5.4);
        
    \end{scope}
    \node[above] at (90:2) {$\bN$};
    \node[below] at (-137:2.1) {$\mathbb{N}$};
    
    \draw (-1.1,0.8) -- (-2.2,0.8) node[left] {$\kbN$};
    \draw (-100:1.6) -- (-90:2.4) node[below] {a thick set};
    \draw (0:1.6) -- (0:2.4) node[right] {a syndetic set};
\end{tikzpicture}
\end{center}
\caption{Illustration of the ultrafilters characterizing syndetic and thick sets. The minimal left ideals are indicated as horizontal lines inside the set $\kbN$. Thick sets contain a whole line while syndetic sets intersect every line.\label{figure-syndetic-thick}}
\end{figure}
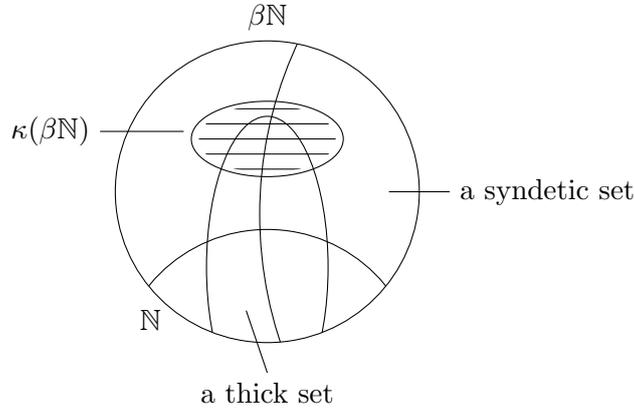

It is also possible to characterize syndetic and thick sets in terms of ultrafilters in the minimal ideal $\kbN$. Of course, since neither of these properties is partition regular, checking whether the closures contain a single ultrafilter is not enough. But we can look for sets of ultrafilters and obtain the following statements, illustrated in Figure \ref{figure-syndetic-thick}.

\begin{proposition}[Ultrafilters characterizing thick sets]
A set $A$ of natural numbers is \emph{thick} if and only if its closure $\close A$ contains \emph{at least one} minimal left ideal $\bN + p \subseteq \kbN$.
\end{proposition}
\begin{proposition}[Ultrafilters characterizing syndetic sets]
A set $A$ of natural numbers is \emph{syndetic} if and only if its closure $\close A$ has \emph{nonempty intersection with every} minimal left ideal $\bN + p \subseteq \kbN$.
\end{proposition}

%: FIGURE: syndetic and thick

We prove these statements and prepare the proof of the theorem by collecting several alternative criteria for being syndetic, thick or piecewise syndetic in the following lemma.

\begin{lemma}[Alternative criteria: thick, syndetic, piecewise syndetic]
% Overfull hbox OK
\label{prop-criteria-thick}
\label{prop-criteria-syndetic}
A set $A$ of natural numbers …
\begin{enumerate}
\item … is \emph{thick} if and only if all finite intersections of shifts are nonempty, \[\bigcap_{n=1}^N (A-n) \neq \emptyset \quad\text{for all } N\in\N.\]
\item … is \emph{thick} if and only if there exists an ultrafilter $p\in\bN$ such that $A-n \in p$ for all shifts $n\in\N$.
\item … is \emph{thick} if and only if its closure $\close A$ contains a minimal left ideal, i.e. $\bN + p \subseteq \close A$ for some ultrafilter $p\in\bN$.
\item … is \emph{syndetic} if and only if its complement $A^c$ is \emph{not thick}.
\item … is \emph{syndetic} if and only if the closure of its complement $\close A^c$ contains no minimal left ideal.
\item … is \emph{syndetic} if and only if its closure $\close A$ has nonempty intersection with every minimal left ideal.
\item … is \emph{piecewise syndetic} if and only if some finite union of shifts $\Union_{n=1}^N (A-n)$ is thick.
\item … is \emph{piecewise syndetic} if and only if there exists an ultrafilter $p\in\bN$ such that the set $A-p$ is syndetic.
\end{enumerate}
\end{lemma}
\begin{proofx}
\begin{enumerate}
\item If the set $A$ is thick, then it contains an interval $[M+1,M+N] \subseteq A$ for each size $N\in\N$. Clearly, the mentioned intersection is nonempty because it contains at least the number $M$.

Conversely, if a set $A$ is not thick and only contains intervals up to a certain size $N$, then taking the intersection of $N+1$ shifts is definitely empty.
\item Consider the infinite intersection of closed sets $B = \bigcap_{n=1}^\infty \close{A-n}$. The previous criterion tells us that every finite intersection is nonempty. But as the space $\bN$ is compact, the infinite intersection must be nonempty as well. In particular, it contains an element $p\in B$, which is the ultrafilter we desire.
\item The previous criterion says that $\N+p \subseteq \close A$. Taking the closure yields $\close{\N+p}=\bN+p \subseteq \close A$. While $\bN+p$ is a left ideal, it is not necessarily minimal, but we know that it certainly \emph{contains} a minimal right ideal.
\item Take complements: $\Union_{n=1}^N (A-n) = \N \iff \bigcap_{n=1}^N (A^c-n) = \emptyset$.
\item Negate criterion 3 and apply it to criterion 4.
\item Another way of writing the previous criterion.
\item Assume that the set is an intersection $A=S \cap T$ of a thick set $T$ and a syndetic set $S$ with $\Union_{n=1}^N (S-n) = \N$. Then, its union of shifts contains the set
\[ \Union_{n=1}^N (A-n) = \Union_{n=1}^N (S-n)\cap(T-n) \supseteq \Union_{n=1}^N \left((S-n)\cap\left(\bigcap_{n=1}^N (T-n)\right)\right) = \bigcap_{n=1}^N (T-n)  .\]
% Overfull hbox OK
But the intersection of shifts of a thick set is again a thick set.

Conversely, let $A$ be a set such that a finite union of its shifts is thick. It is easy to find a slightly larger thick set $T$ such that $\bigcap_{n=1}^N (T-n) \subseteq \Union_{n=1}^N (A-n)$, for instance by enlarging the intervals a bit. Then, the set $S := A \cup (\N \setminus T)$ fulfills the condition $A = S \cap T$ and it is also syndetic because
\[ \Union_{n=1}^N (S-n) \subseteq \Union_{n=1}^N (A-n) \cup \Union_{n=1}^N (\N\setminus T-n) = \Union_{n=1}^N (A-n) \cup 
\left(\N \setminus \bigcap_{n=1}^N (T-n) \right) = \N .\]
\item Combine the previous criterion and criterion 2 to obtain that a set $A$ is piecewise syndetic if and only if
\[ \Union_{n=1}^N (A-n) - p = \N .\]
But we can interchange the order of shifts $(A-n)-p = (A-p)-n$ and obtain that $A-p$ is syndetic.\hfill\qedsymbol
\end{enumerate}
\end{proofx}

With these preliminaries in place, the following lemma gives the key step to the characterization of piecewise syndetic sets.

\begin{lemma}[Shift by element from minimal ideal is syndetic]
\label{prop-shift-minimal-ultrafilter-syndetic}
Let $A$ be a set of natural numbers and assume that its closure $\close A$ contains an ultrafilter $p\in\kbN$ from the minimal ideal. Then, the set $A-p$ is syndetic. 
\end{lemma}
\begin{proofx}
Since the ultrafilter $p$ is contained in a minimal left ideal, we know that for every ultrafilter $q\in\N$, there exists an ultrafilter $r_q\in\bN$ such that
\[  r_q+q+p = p .\]
Since $p\in\close A$, this implies
\[ r_q+q+p \in \close A \iff A \in r_q+q+p \iff (A-p)-q \in r_q .\]
By the permanence principle, there exists a natural number $n_q \in \N$ such that
\[ (A-p)-q \in n_q \iff (A-p)-n_q \in q\]
Taking the union over all ultrafilters $q\in\bN$, we can cover the space $\bN$ with countably many open sets
\[ \Union_{n\in\N} \close{(A-p)-n} = \bN .\]
But since the space $\bN$ is compact, it can already be covered by finitely many such open sets, and we obtain that $A-p$ is syndetic,
\[ \Union_{n=1}^N ((A-p)-n) = \N .\]
\qedEquationHack{0.98cm}
\end{proofx}

We can now complete the proof of the theorem.

\begin{proof}[Ultrafilters characterizing piecewise syndetic sets]

\Backward:
Consider a basic open set $\close A$ and assume that it contains an ultrafilter $p\in\kbN$. Lemma \ref{prop-shift-minimal-ultrafilter-syndetic} tells us that the ultrafilter shift $A-p$ is syndetic and Lemma \ref{prop-criteria-syndetic} tells us that this means that $A$ is piecewise syndetic.

\Forward:
Assume that the set $A$ is piecewise syndetic, which means that we can write it as an intersection $A=T\cap S$ of a thick set $T$ and a syndetic set $S$.

Lemma \ref{prop-criteria-syndetic} about alternative criteria for being thick or syndetic tells us that the closure $\close T$ contains a minimal ideal, while the closure $\close S$ has nonempty intersection with every minimal ideal. In particular, the closure $\close S$ intersects the minimal ideal contained in the set $\close T$. Hence, there exists an ultrafilter $p\in\kbN$ contained in the intersection $\close S \cap \close T = \close A$.
\end{proof}

%%%%%%%%%%%%%%%%%%%%%%%%%%%%%%%%%%%%%%%%%%%%%%%%%%%%%%%%%%%%%%%%%%%%%%
\section{Counting measures on \titlebN}
\label{chapter-measure}

To interpret Szemerédi's theorem in terms of ultrafilters, we need a notion of density on the space of ultrafilters $\bN$, which will be given by so-called ``\emph{counting measures}''. One advantage of the space of ultrafilters over the set of natural numbers $\N$ is that we no longer have to content ourselves with a subadditive set function like the \density{}, instead we can now construct proper $\sigma$-additive \emph{measures}. 

If you are familiar with Haar measures for compact topological groups, you may notice a certain similarity, but keep in mind that the compact space $\bN$ has neither a properly continuous semigroup operation, nor are counting measures in any way unique.

%%%%%%%%%%%%%%%%%%%%%%%%%%%%%%%%%%%
\subsection{Construction}
\label{section-measure-construction}

The problem with defining the density of a set $A\subseteq\N$ was that, in general, the sequence of fractions
\[ \alpha_N = \frac{|A\cap[1,N]|}{N} = \frac1N\sum_{k=1}^N \1_A(k) \in[0,1]\]
does not converge for $N\to\infty$. For instance, consider the set $A$ consisting of all natural numbers with an odd number of decimal digits: its sequence of fractions will eventually oscillate between $1/11$ and $10/11$. When defining the \density{}, we ``solved'' this problem by taking the limes superior, but this comes at a terrible price: the \density{} is no longer additive, only subadditive.

Fortunately, ultrafilter limits can help us to restore additivity. We simply choose an ultrafilter $p$ and declare the ultrafilter limit $\lim_{N\to p} \alpha_N$ to denote the density of the set $A$. Unlike the limes superior, ultrafilter limits are compatible with addition (Proposition \ref{ultrafilter-limit-addition}) and the resulting notion of density will be properly \emph{additive}. Of course, the price we now pay is that there is no unique choice for the ultrafilter $p$.

We cannot expect such a density to be \emph{$\sigma$-additive} and a proper \emph{measure}, however, because the space of natural numbers $\N$ is countable and thus wholly unsuitable for playing the role of a measure space. But we have seen that the space of ultrafilters $\bN$ is more appropriate in the context of Ramsey theory anyway (Section \ref{section-ramsey-ultrafilter} ), so it is only natural to apply the notion of density to sets of ultrafilters. Since $\bN$ is a compact topological space, we will obtain proper measures, called ``\emph{counting measures}'', which we will now define.

Similarly to the upper Banach density, we allow the fractions to be taken over more general intervals than just $[1,N]$.

\begin{definition}[Counting measure on $\bN$]
\label{definition-counting-measure}
A measure $\mu$ on $\bN$ is called a \define{counting measure} if it is a Radon measure that is given by an ultrafilter limit
\[  \mu(\close A) = \lim_{(N,M)\to q} \frac{|A\cap[M+1,M+N]|}{N}  \quad\text{for all }A\subseteq\N ,\]
where the ultrafilter $q\in\beta(\N^+\times\N)$ fulfills two conditions:
\begin{properties}
\item The interval length $N$ tends to infinity, i.e.
\[ \lim_{(N,M)\to q} \frac1N = 0 .\]
% toto: layout
\item The purely technical requirement that either the interval starting points are bounded $\lim_{(N,M)\to q} M = M_0 < \infty$, or there exists a set $I=\{(N_1,M_1),$ $(N_2,M_2) \dots\}$ of \emph{disjoint} $M_1+N_1 < M_2, \dots$ and \emph{growing} $N_1 < N_2 < \dots$ intervals such that the ultrafilter $q$ is contained in the closure $\close I$.
\end{properties}
\end{definition}

Being a Radon measure is a regularity condition that allows us to calculate the measures of arbitrary measurable sets in terms of the basic open sets, for which we have given a more or less explicit formula. In particular, we obtain a formula for the measure of a closed set.

\begin{proposition}[Measure of a filter]
\index{filter!measure of a filter}
\label{prop-measure-filter}
Let $F$ be a closed subset of $\bN$, corresponding to a filter $\F = \{A\subseteq\N : F \subseteq\close A \}$. Then, the measure of $F$ can be calculated as
\[  \mu(F) = \inf_{A\in\F} \mu(\close A) \]
\end{proposition}
\begin{proof}
Since $\mu$ is a Radon measure, the compact set $F$ fulfills
\[  \mu(F) = \inf_{U \text{ open},\ F \subseteq U} \mu(U) .\]
See \textcite{Elstrodt:2005p567} Chapter VIII, Corollary 1.2, page 314. But every open set is a union of basic open sets and it is sufficient to take the infimum over all basic open sets as indicated.
\end{proof}

Given an ultrafilter $q$, one could probably show that the corresponding counting measure is well-defined and uniquely determined. However, it is easier to construct these measures from positive linear functionals. The following proposition will supply us with all the counting measures we need.

\begin{proposition}[Counting measure from \density{}]
\index{upper Banach density}
Let $A$ be a set of natural numbers with \density{} $\d(A)$. Then, there exists a counting measure $\mu$ on $\bN$ with the property
\[  \mu(\close A) = \d(A)  .\]
\end{proposition}
\begin{proof}
First, we construct the ultrafilter $q$ used in the definition of counting measures.

It is clear from the definition of \density{} that there exists a sequence of intervals $I_j = [M_j+1,M_j+N_j]$ of unbounded length $N_j\to\infty$ whose sequence of ratios ${|A\cap I_j|}/{N_j}$ converges to $\d(A)$. Now, apply a variant of Lemma \ref{prop-subsequences} about ultrafilters and subsequences to see that any non-principal ultrafilter $q\in\close{\{(N_1,M_1),(N_2,M_2),\dots\}} \subseteq \beta(\N^+\times\N)$ will give
\[  \lim_{(N,M)\to q} \left(\frac{|A\cap [M+1,M+N]|}N, \frac 1N\right) = (\d(A),0) .\]
It is also clear that we can fulfill the second requirement on the ultrafilter $q$ by first passing to a subsequence of intervals if necessary.

Now, use such an ultrafilter $q$ to define a positive linear functional $I : \boundedSequences \to \R$ on the space of bounded sequences as follows
\[ I(f) := \lim_{(N,M)\to q} \frac 1N\sum_{k\in [M+1,M+N]} f(k).\]
The ultrafilter limit exists because the average of a bounded sequence $f(k)$ is again bounded; clearly we have $|I(f)| \leq ||f||_\infty$. Linearity and positivity are also obvious. Note that characteristic functions $\1_B$ of sets $B\subseteq\N$ are also bounded sequences and that we have $I(\1_A) = \d(A)$.

Remembering the correspondence $\boundedSequences = \Continuous(\bN)$ from Proposition \ref{prop-bounded-sequences}, we can apply \emph{Riesz' representation theorem} (\textcite{Elstrodt:2005p567}, chapter VIII, theorem 2.5, page 335) to this linear functional and obtain a Radon measure $\mu$ on the space of ultrafilters with the property $\mu(\close B) = I(\1_B)$. This makes it clear that $\mu$ is a counting measure and that we have $\mu(\close A) = I(\1_A) = \d(A)$ as desired.
\end{proof}

As you would expect, counting measures are invariant under finite shifts.
\begin{proposition}[Shift invariance]
\index{counting measure!shift invariance}
Counting measures are shift invariant, i.e.\ we have
\[  \mu(F) = \mu(F-1)  \quad\text{where } F-1:=\{p\in\bN : p+1\in F \}\]
for every measurable set of ultrafilters $F\subseteq\bN$.
\end{proposition}
\begin{proof}
This is very much a consequence of the fact that the interval length tends to infinity. After all, for a basic open set $F=\close B$ with $B\subseteq\N$, we have
\[ \left|\frac{|(B-1)\cap[M+1,M+N]|}{N} - \frac{|B\cap[M+1,M+N]|}{N}\right| \leq \frac2 N\]
and taking the ultrafilter limit $(N,M)\to q$ gives $|\mu(\close B) - \mu(\close B-1)| = 0$ as desired.

To see that this equality extends to all measurable sets, note that we can define a new measure $\mu'(F) := \mu(F-1)$ since shifting commutes with all $\sigma$-algebra operations. In other words, $\mu'$ and $\mu$ are Radon measures that coincide on all basic open sets; they must be equal everywhere.
\end{proof}

To summarize, we have translated the notion of \density{} of a set of natural numbers into a bona-fide measure on the space of ultrafilters $\bN$.

%%%%%%%%%%%%%%%%%%%%%%%%%%%%%%%%%%%
\subsection{Application: Sets of differences}
Let us demonstrate the utility of the counting measures on $\bN$ by proving the following theorem of Jin \parencite{Jin:2002p582}.

\begin{theorem}[Jin. Set of differences is piecewise syndetic]
\index{Jin's theorem}
Let $A$ and $B$ be two sets of natural numbers with positive upper Banach density. Then, the set of differences\footnote{Actually, Jin has proved this for the set of sums $A+B$, but that is of negligible importance here.}
\[A-B := \{a-b : a\in A,b\in B\}\cap\N\]
is \emph{piecewise syndetic}.
\end{theorem}

Jin's theorem provides a link between two notions of largeness: a subset of the natural numbers can be large in the sense that it has positive density, or it can be large in the sense that it is piecewise syndetic. As we will prove later (Proposition \ref{prop-density-not-piecewise-syndetic}), these two notions are not the same: a set can have high density while not being piecewise syndetic; but the theorem demonstrates that they are not completely independent.

The proof presented here is due to \textcite{Beiglbock:2009p581}. It proceeds by reducing the theorem to the following easy special case.

%:	C-C syndetic

\begin{proposition}[Set of self-differences is syndetic.]
Let $C$ be a set of natural numbers with positive upper Banach density. Then, the set of ``self''-differences $C-C$ is syndetic.
\end{proposition}
\begin{proof}
Consider the shifted sets $C-n$ for $n\in\N$. The idea is that these sets intersect each other very often. More precisely, we claim that we can choose a collection of shifts $\{\vector nk\}$ whose size $k$ is \emph{maximal} among those collections who respect the condition that the shifted sets $(C-n_1),(C-n_2),\dots,(C-n_k)$ should be \emph{mutually disjoint}.

To see this, choose a counting measure $\mu$ such that $\mu(\close C)=\d(C) > 0$. By finite additivity and shift invariance, we have
\[  \mu(\close{C-n_1} \uplus \close{C-n_2} \uplus \dots \uplus \close{C-n_k}) = k\cdot\mu(\close C) \leq \mu(\bN)=1  .\]
Hence, the size of any collection of mutually disjoint shifted sets is bounded by $k \leq 1/\mu(\close C)$.

Now, assume indeed that the shifts $C-n_i$ form a maximal collection of mutually disjoint sets. This implies that every other shifted set $C-m$ must intersect one of the $C-n_i$, i.e.\ $(C-m)\cap(C-n_i) \neq \emptyset$ for one shift $n_i$. Another way of writing this is $m \in (C-C)+n_i$. Taking the union over all $m$ gives
\[  \bigcup_{i=1}^k ((C-C) + n_k) = \N  ,\]
which means that $C-C$ is syndetic.
\end{proof}

%:	ultrafilter shift

The key idea for proving the general case is to find a shift $n$ such that the intersection $C=(A-n)\cap B$ also has positive upper Banach density. Then, the theorem follows by applying the special case to the set $C$. Unfortunately, this will not work if we only consider natural numbers $n\in\N$, but fortunately we can allow the shift to become an \emph{ultrafilter} $n\in\bN$.

\begin{lemma}[Intersection with shift by an ultrafilter]
Let $A,B$ be two sets of natural numbers. Then, there exists an ultrafilter $p\in\bN$ with
\[  \d((A-p)\cap B) \geq \d(A)\cdot\d(B)  .\]
\end{lemma}
\begin{proofx}
\newcommand{\intIn}{\frac1{|I_n|}\sum_{k\in I_n}}%
\newcommand{\intbN}[1]{\int_\bN #1 \ d\mu}%
Choose a counting measure $\mu$ on the space of ultrafilters such that $\mu(\close A)=\d(A)$. Furthermore, let $I_n \subseteq\N$ be a sequence of intervals of increasing length that realizes the supremum in the definition of the upper Banach density of $B$, i.e.\ 
$  \d(B) = \lim_{n\to\infty} \frac{|B\cap I_n|}{|I_n|}  .$
Additionally, let $\1_X$ denote the characteristic function of a closure $\close X\subseteq\bN$. For instance, we can write
$  \intbN{\1_A(p)} = \mu(\close A)$ and $\lim_{n\to\infty}\intIn \1_{B}(k) = \d(B)  .$

%\to do{? Steinhaus lemma. Appearance of convolution is not surprising.} This remark is out.
Now, consider the following integral, which resembles the integral of a \emph{convolution} of the functions $\1_A$ and $\1_B$:
\[  \intbN{\intIn \1_A(k+p)\1_B(k)}  .\]
Since all functions are bounded and the space has finite measure, we can apply Fatou's lemma in the following form:
\[  \intbN{ \limsup_{n\to\infty} \intIn \1_A(k+p)\1_B(k)}
    \geq  \limsup_{n\to\infty}  \intbN{\intIn \1_A(k+p)\1_B(k)}  .\]

For the left-hand side, we write the product of characteristic functions as
\[  \1_A(k+p)\1_B(k) = \1_{A-p}(k)\1_B(k) = \1_{(A-p)\cap B}(k)  \]
and reason that for each $p$, the limes superior over these particular intervals $I_n$ must be smaller than the upper Banach density
\[  \intbN{ \d((A-p)\cap B) } \geq \intbN{ \limsup_{n\to\infty} \intIn \1_{(A-p)\cap B}(k) } = \text{left-hand side}.\]

For the right-hand side, we write the product as
\[  \1_A(k+p)\1_B(k) = \1_A(p+k)\1_B(k) = \1_{A-k}(p)\1_B(k)  \]
and integrate over $p$, noting that $\mu$ is shift-invariant
\[\begin{split}
  \text{right-hand side} &= \limsup_{n\to\infty} \intIn\1_B(k) \intbN{\1_{A-k}(p)} \\
  &= \limsup_{n\to\infty}\intIn\1_B(k) \mu(\close A) = \d(B)\mu(\close A)
.\end{split}\]

Putting both together and remembering that $\mu(\close A)=d(A)$, we obtain the inequality
\[  \intbN{ \d((A-p)\cap B) } \geq \d(A)\cdot\d(B) .\]
Since $\mu(\bN)=1$, this is only possible if there exists at least one ultrafilter $p$ for which
\[  \d((A-p)\cap B) \geq \d(A)\cdot\d(B) .\]
\qedEquationHack{0.98cm}
\end{proofx}

%:	A-B syndetic

\begin{proof}[Jin's theorem]
By the lemma, there is an ultrafilter $p$ such that the set $C=(A-p)\cap B$ has positive upper Banach density. Hence, the set of differences $C-C$ is syndetic.

Now, consider the ultrafilter shift $(A-B)-p$. We have
\[ (A-B)-p \supseteq (A-p) - B \supseteq C - C \]
which implies that it is also syndetic. The first inclusion is justified by
\[ \begin{split}
  n \in (A-p)-B &\iff \exists b\in B.\ n \in (A-p) - b \iff \exists b\in B.\ A-b \in n+p \\
  &\implies A-B \in n+p \iff n \in (A-B) - p
.\end{split}\]
But Lemma \ref{prop-criteria-syndetic} about criteria for being piecewise syndetic tells us that the set $A-B$ is piecewise syndetic if there exists an ultrafilter such that the shift $(A-B)-p$ is syndetic. This completes the proof.
\end{proof}

%%%%%%%%%%%%%%%%%%%%%%%%%%%%%%%%%%%
\subsection{Measuring \aprich\ ultrafilters}
\label{section-szemeredi-ultrafilter}

We now present the link between Szemerédi's theorem (\ref{theorem-szemeredi}) and the size of the set of \aprich\ ultrafilters $AP_k$ (\ref{definition-aprich}).

The key observation is that the property of having positive measure is \emph{partition regular}. Namely, any partition $A = \partition{A}$ implies that
\[  \mu(\close A) = \mu(\close A_1) + \mu(\close A_2) + \dots + \mu(\close A_r)  \]
and if the closure $\close A$ has positive measure, then one of the parts $\close A_i$ must have positive measure, too. By the ultrafilter construction lemma, this means that any set with positive measure must contain an ultrafilter $p$ whose basic open neighborhoods all have positive measure.

But now, Szemerédi's theorem tells us that each of these open neighborhoods contain an arithmetic progression of length $k$. This means that the ultrafilter $p$ is also \aprich, $p\in AP_k$. In other words, any set with positive measure not only contains a single arithmetic progression, but a whole \emph{ultrafilter} $p$ that is ``full of arithmetic progressions''.

Let us recast this observation above in terms of the following standard notion.
\begin{definition}[Support of a measure (on $\bN$)]
The \define{support of a measure} $\mu$ on $\bN$ is defined to be the set of points whose basic open neighborhoods all have positive measure:
\[  \supp\mu = \{ p \in\bN : \mu(\close A) > 0 \text{ for all } p\owns A \}  .\]
\end{definition}
Thus, we have just argued that any set with positive measure intersects the support of $\mu$, and that Szemerédi's theorem implies that $\supp\mu \subseteq AP_k$. This looks like a strong lower bound on the the size of $AP_k$, because we expect $\supp\mu$ to be large. The following proposition tells us how large.

\begin{proposition}[Size of the support of a measure]
For any Radon measure $\mu$ on $\bN$, we have $\mu(\supp\mu) = \mu(\bN)$.
\end{proposition}
\begin{proofx}
Put differently, we expect that the complement $(\supp\mu)^c$ is a null set. But we have argued above that every open set $\close A$ with positive measure already intersects $\supp\mu$, so every $\close A$ that is contained in $(\supp\mu)^c$ must be a null set. The inner regularity of our Radon measure now implies that
\[  \mu((\supp\mu)^c) = \sup_{\close A\subseteq (\supp\mu)^c} \mu(\close A) = 0  .\]
\qedEquationHack{0.98cm}
\end{proofx}

Taken together, these observations yield the following result.

%:	Szemerédi

\begin{theorem}[Szemerédi's theorem is equivalent to $AP_k$ being large]
The following are equivalent:
\begin{itemize}
\item Szemerédi's theorem\index{Szemeredi's theorem@Szemerédi's theorem}: If $\d(A) > 0$, then $A$ contains an arithmetic progression of length $k$.
\item The set of \aprich\ ultrafilters has full measure: $\mu(AP_k) = \mu(\bN)$ for any counting measure $\mu$ on $\bN$.
\end{itemize}
% Overfull hbox OK
\end{theorem}
\begin{proof}
\Forward :
Assuming Szemerédi's theorem, we have already argued that $\supp\mu \subseteq AP_k$. But since $\mu(\supp \mu) = \mu(\bN)$, this implies $\mu(AP_k) = \mu(\bN)$ as desired.

\Backward :
Let $A\subseteq\N$ be a set with positive upper Banach density. We choose a counting measure on the space of ultrafilters such that $\mu(\close A)=\d(A)>0$.

By assumption, $\mu(AP_k)=\mu(\bN)$, so its complement $AP_k^c$ must be a null set. But of course, this makes it impossible for a set $\close A$ with nonzero measure to be fully contained in the complement. In other words, the closure $\close A$ must intersect $AP_k$ and hence, the set $A$ contains an arithmetic progression of length $k$.
\end{proof}

So, unlike van der Waerden's theorem, for which the existence of a single ultrafilter $p\in AP_k$ was sufficient, any attempt to prove Szemerédi's theorem requires us to ``cough up`` a large amount of \aprich\ ultrafilters. That being said, we did identify more than one \aprich\ ultrafilter; namely, we showed that the minimal two-sided ideal $\kbN$ is contained in the set $AP_k$ in Section \ref{kappa-is-aprich}. Might that be enough to imply Szemerédi's theorem? No, as the following result shows.

%:	k(bN) is a null set

\begin{theorem}[The minimal two-sided ideal $\kbN$ is a null set]
\index{ideal!minimal two-sided ideal@minimal two-sided ideal $\kbN$}
The smallest two-sided ideal $\kbN$ in $\bN$ is a null set, i.e.\ $\mu(\kbN) = 0$ for any counting measure $\mu$.
\end{theorem}

In particular, showing that $\kbN \subseteq AP_k$ is not enough to prove Szemerédi's theorem, because that only implies the trivial lower bound $0 = \mu(\kbN) \leq \mu(AP_k)$.

\begin{proofx}
The subsequent proposition will construct sets $A$ that are not piecewise syndetic, but have arbitrarily high density $<1$. Since a set is piecewise syndetic if and only if it intersects $\kbN$, these sets $\close A$ will be contained in the complement $\kbN^c$. Since they have arbitrarily high density, we conclude
\[   \mu(\kbN^c) = \sup_{\close A \subseteq \kbN^c} \mu(\close A) \geq 1 = \mu(\bN)  .\]
\qedEquationHack{0.98cm}
\end{proofx}

\begin{proposition}[Sets of high density that are not piecewise syndetic]
\index{piecewise syndetic}
\index{upper Banach density}
\label{prop-density-not-piecewise-syndetic}
For every desired density $\alpha < 1$, there exists a set $A$ with counting measure $\mu(\close A) \geq \alpha$ which is not piecewise syndetic.
\end{proposition}
\begin{proof}
First of all, we need to understand what it means for a set $A$ of natural numbers to not be piecewise syndetic. It means that no union of $\Union_{k=1}^d (A-k)$ may be thick, which implies that all intervals fully contained in this union have length no larger than some number $l(d)\in\N$. In other words, we have the following equivalence:
\begin{itemize}
\item The set $A$ is not piecewise syndetic.
\item For every gap width $d\in\N$, there exists an interval length $l(d)$ such that every intersection of the set $A$ with an interval $[M+1,M+l(d)] \subseteq \N$ contains a gap of length greater or equal than $d$. An illustration is given in Figure \ref{figure-not-piecewise-syndetic}.
\end{itemize}

%: FIGURE: not piecewise syndetic set
\begin{figure}[ht]
\begin{center}
\begin{tikzpicture}[yscale=0.5,xscale=0.66,scale=0.6]
	\def\block#1#2{\draw (#1,-1) rectangle +(#2,2)}
	\def\number#1[#2]{\node[#2] at (#1+0.5,0) {$#1$}}

	%\draw[help lines] (0,-2) grid (10,2);
	\node[left] at (0.5,0) {$A = \Big\{$};
	\block{1}{3};
	\number{1}[black];
	\number{2}[black!70];
	\number{3}[black!60];
	
	\node at (5.6,0) {$…$};
	
	% interval
	\block{7}{4};
	\block{12}{5};
	\block{20}{3};
	\block{24}{4};
		
	\def\gapsize#1#2{
	\draw (#1,-2) -- +(0,-0.5) -- node[midway,below,style={font=\footnotesize}] {some gap has size $\geq d$} (#2,-2.5) -- +(0,0.5)}
	\gapsize{17}{20};

	\def\intervals#1#2#3{
	\draw (#1,-4.7) -- +(0,-0.5) -- node[midway,below,style={font=\footnotesize}] {#3} (#2,-4.7-0.5) -- +(0,0.5)}
	\intervals{9}{30}{when looking at any interval of size $\geq l(d)$};
	
	\path (30.6,0) node {$\dots\dots$} +(2.5,0) node {$\Big\}$};
\end{tikzpicture}
\end{center}
\caption{A set $A$ is not piecewise syndetic if and only if a gap of length at least $d$ can be found in each interval that has size at least $l(d)$.\label{figure-not-piecewise-syndetic}}
\end{figure}
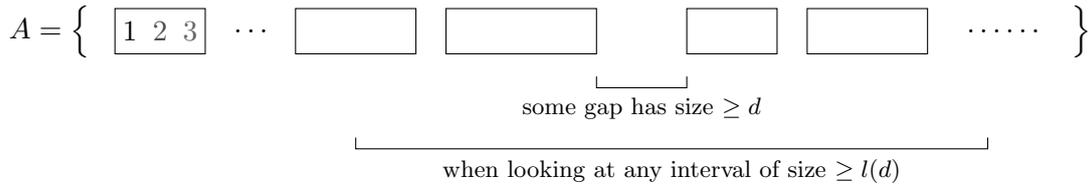

We now want to construct a set $A$ with both this property and the property that $|A \cap [1,N]|/N \geq \alpha$ for all interval lengths $N\in\N^+$. For simplicity, we begin with $\alpha = 1/2$ and later argue that we can also achieve higher densities.

Note that counting measures can also involve intervals other than $[1,N]$, so this set $A$ may not be the one we actually want. But this is a harmless problem that we will deal with at the very end.

To construct the set $A$, a ``self-similar'' or ``fractal'' process seems to be most suitable, since fractals tend to have large gaps while also mainting a high density. Hence, consider the following sequence of words $A_n$, made from two symbols ``$0$'' and ``$1$'':
\begin{align*}
 A_0     &= 1 \\
 A_1     &= 110 \\
 A_2     &= 1101100 \\
 A_3     &= 110110011011000 \\
 &\vdots \\
 A_{n+1} &= A_n A_n 0	\quad\text{ for all } n\geq 1
\end{align*}
Since each word is a prefix of the subsequent one, the set
\[ \begin{split}
 A = \{ k\in\N : &\text{ the $k$-th symbol in the word $A_n$ is equal to ``$1$''}
\\& \text{ for sufficiently large $n\in\N$} \}
\end{split} \]
is well-defined. In other words, we interpret the set $A\subseteq\N$ as an infinite sequence $A\in\{0,1\}^\N$ of the symbols $0$ and $1$ and construct it by a series of finite approximations $A_n$.

\def\length{\mathop{\textit{length}}}
First, let us argue that the set $A$ is not piecewise syndetic. To see that, note that the each words $A_n$ ends with exactly $n$ zeroes, i.e.\ it has the form
\[ A_n = B_n\underbrace{000…0}_\text{$n$ zeroes} \]
for some word $B_n$ that may contain the symbol ``$1$''. This is immediately clear from induction. But the recursive nature makes it also clear that the set $A$ consists of blocks made from the word $A_n$:
\def\B{B_n\underbrace{000…0}_n}
\[ A = A_nA_n0A_nA_n00… = \B\B0\B\B00…\]
This means that the set $A$ consists of blocks made from the word $B_n$ interspersed by blocks consisting of $n$ or more zeroes. Hence, any interval of length $l(n) \geq \length(B_n)+2n$ will contain a gap of length at least $n$. This proves that the set $A$ is not piecewise syndetic.

\def\ones{\mathop{\textit{ones}}}
Now, let us argue that the set $A$ fulfills the density property $|A\cap[1,N]|/N \geq 1/2$. To see that, first note the following properties
\begin{itemize}
\item The length $\length(A_n)$ of the word $A_n$ is equal to
\[ \length(A_n) = 2^{n+1} - 1.\]
After all, we have the recursive equations $\length(A_0)=1$ and $\length(A_{n+1})=2\cdot\length(A_n)+1$.
\item The number of times $\ones(A_n)$ that the symbol ``$1$'' appears in the word $A_n$ is equal to
\[ \ones(A_n) = 2^n\]
because of the recursive equations $\ones(A_0)=1$ and $\ones(A_{n+1})=2\cdot\ones(A_n)$.
\end{itemize}
Hence, if the interval length $N$ is exactly the length $N = \length(A_n)$ of a word, we have
\[ \frac{|A\cap[1,N]|}{N} = \frac{\ones(A_n)}{\length(A_n)} = \frac{2^n}{2^{n+1}-1} \geq \frac 12 .\]
But if the interval length $N$ is not exactly the length of a word, we can use the relation $A_{n+1} = A_nA_n0$ to split the interval into two parts and use induction, like this: assume that for all interval sizes $N \leq \length(A_n)$, we have shown that the density is $\geq 1/2$. Now, consider a size $\length(A_n) < N < \length(A_{n+1})$. We can use the induction hypothesis and the self-similar shape to conclude that
\begin{align*}
 \frac{|A \cap [1,N]|}{N}
  &= \frac1N\Big(|A\cap[1,\length(A_n)]|+|A\cap[\length(A_n)+1,N]| \Big) \\
  &= \frac1N\Big(|A\cap[1,\length(A_n)]|+|A\cap[1,N-\length(A_n)]| \Big) \\
  &\geq \frac1N\left( \frac12\length(A_n) + \frac12(N-\length(A_n)) \right)
  = \frac12
.\end{align*}
This concludes the argument about the density.

To obtain densities $\alpha$ greater than $1/2$, we have to start with the word $A_0=11\dots1$ consisting of $k$ ones instead of just a single symbol ``$1$''. The argument showing that the set $A$ is not piecewise syndetic is left unchanged while the calculation of the density gives
\begin{align*}
\length(A_n) &= 2^n(k+1)-1\\
\ones(A_n)  &= 2^nk \\
\frac{|A \cap [1,N]|}{N} &\geq \frac{\ones(A_n)}{\length(A_n)} = \frac{k}{(k+1)-\frac1{2^n}}\geq \frac{k}{k+1}
\end{align*}
which can be chosen arbitrarily close to $1$.

Finally, we have to discuss what the density property $|A\cap [1,N]|/N \geq \alpha$ means for the measure $\mu(\close A)$. Clearly, the counting measure depends on the ultrafilter $q$ associated to it. To obtain a lower bound for the measure $\mu(\close A)$, we consider a collection of intervals $I = \{(N_1,M_1),(N_2,M_2),\dots\}$ with $q \in \close I$ and simply show that the set $A$ has density at least $\alpha$ on all of these intervals, $|A\cap [M_i+1,M_i+N_i]| / N_i \geq \alpha$. Then, the same inequality must hold for the ultrafilter limit, see also Lemma \ref{prop-limit-restriction}.

To find a good collection of intervals $I$, we will make use of the second requirement on the ultrafilter $q$ in the definition of counting measures \ref{definition-counting-measure}.

The first case is that the interval starting points are bounded, $\lim_{(N,M)\to q} M = M_0 < \infty$. The definition of the limit implies that $q \in \close{\{(N,M_0) : N\in \N^+\}}$. In other words, the counting measure $\mu(\close A)$ only depends on the densities $|A \cap [M_0+1,M_0+N]|/N$. But when the interval length $N$ becomes much larger than the fixed starting point $M_0$, the difference of densities
\[ \left| \frac{|A \cap [M_0+1,M_0+N]|}N - \frac{|A \cap [1,N]|}N \right| \leq 2\frac{M_0}N \]
becomes negligible and we conclude $\mu(\close A) \geq \alpha$.

In the other case, where the interval starting points are not bounded, the measure $\mu(\close A)$ may, in fact, vanish because the set $A$ contains arbitrarily large gaps and the intervals from the collection $I$ might fall exactly into these gaps. We cannot expect that one and the same set $A$ works for all counting measures!

We have to build a new set $A'$ that is not piecewise syndetic and fulfills $\mu(A') \geq \alpha$. Fortunately, we have required that the ultrafilter $q$ is contained in the closure $\close I$ of a sequence of \emph{disjoint and growing} intervals. We simply fill each of these intervals with a shift of the set $A_n$ surrounded by two gaps of size $n$, i.e.\
\[ A' = \Union_{i\in\N} \left([M_i+1+n_i,M_i+N_i-n_i] \cap (A_{n_i}+M_i+n_i) \right) .\]
Each gap size $n_i$ is chosen as the smallest integer such that $\length(A_{n(i)}) \geq N_i$. Because of $n/\length(A_n) \to 0$, these additional gaps have no impact on density and the set $A'$ will satisfy the density requirement. Furthermore, it is not difficult to check that this set is again not piecewise syndetic because, except for some parts at the beginning, the set will consist of blocks of $A_n$, respectively $B_n$, interspersed by gaps of size at least $n$. This completes the proof.
\end{proof}

%%%%%%%%%%%%%%%%%%%%%%%%%%%%%%%%%%%%%%%%%%%%%%%%%%%%%%%%%%%%%%%%%%%%%%
% Bibliography
\clearpage
\phantomsection
\addcontentsline{toc}{section}{References}
\printbibliography
% Overfull hbox OK

% Notation
%\phantomsection
%\addcontentsline{toc}{section}{Notation}
\clearpage
\phantomsection
\addcontentsline{toc}{section}{List of Symbols}
\section*{List of Symbols}

\begin{tabular}{lp{12.5cm}}
\\ $2^X$ & Powerset, set of subsets of the set $X$.
\\ $A^c$ & Complement of the set $A$ inside some ambient set $X$, $A^c = X\setminus A$. The context should make it clear which set $X$ is meant.
\\ $A \uplus B$ & Union of two disjoint sets $A$ and $B$, i.e.\ $A \uplus B = A \cup B$ with the understanding that $A \cap B = \emptyset$.
\\ $[M,N]$ & Interval of natural numbers; $[M,N]=\{M,M+1,\dots,N-1,N\}$.
\\ $\d(A)$ & Upper Banach density of the set $A\subseteq\N$, see definition \ref{definition-upper-Banach-density}.
\\ $\bN$ & Space of ultrafilters; \StoneCech\ of the set $\N$, see definition \ref{definition-stone-cech}.
\\ $\kbN$ & Minimal two-sided ideal of the left-topological semigroup $\bN$, see definition \ref{definition-minimal-ideal-kappa}.
\\ $AP_k$ & Set of \aprich\ ultrafilters featuring length $k$, see definition \ref{definition-aprich}.
\\ $\close A$ & Basic open set; closure of a set $A\subseteq\N$ in the \StoneCech\ $\bN$, see definition \ref{definition-closure}.
\\ $\closure X$ & Closure of the set $X$ in a topological space.
\\ $\lim\limits_{n\to p} f(n)$ & Limit along the ultrafilter $p$, see definition \ref{definition-ultrafilter-limit}.
\\ $A-p$ & Ultrafilter shift, see definition \ref{definition-ultrafilter-shift}.
\\ $A+B$ & Sumset, $A+B = \{ a+b : a\in A, b\in B\}$.
\end{tabular}

% Index
\clearpage
\phantomsection
\addcontentsline{toc}{section}{Index}

\printindex

% Erklärung
\newpage
\phantomsection
\addcontentsline{toc}{section}{Erkl\"arung}
\section*{Erklärung}
Ich versichere, dass ich die vorliegende Arbeit selbständig und nur unter Verwendung der angegebenen Quellen und Hilfsmittel angefertigt habe,
insbesondere sind wörtliche oder sinngemäße Zitate als solche gekennzeichnet. Mir ist bekannt, dass Zuwiderhandlung auch nachträglich zur Aberkennung
des Abschlusses führen kann.

\vspace{2cm}
\begin{flushright}
Heinrich-Gregor Zirnstein
\end{flushright}

\end{document}